\documentclass[12pt]{article} 
\usepackage{latexsym}
\usepackage[all]{xy}
\usepackage{graphicx}
\def\separation{\medskip}
\def\Z{Z\!\!\!Z}

\def\Q{{\bf Q}}
\def\P{I\!\!P}

\def\k{{\bf k}}

\def\W{{\cal W}}
\def\L{{\cal L}}

\def\I{{\cal I}}

\def\T{{\rm T}}

\def\M{{\cal M}}

\def\?{{\bf ??}}
\def\Proj{{\rm Proj}}

 \def\dim{{\rm dim}}

\def\Im{{\rm Im}}

 \newtheorem{theorem}{Theorem}[section]
\newtheorem{lemma}[theorem]{Lemma} 
\newtheorem{prop}[theorem]{Proposition} 
\newtheorem{definition}[theorem]{Definition} 
\newtheorem{corollary}[theorem]{Corollary}
 
\newtheorem{remark}[theorem]{Remark}

 \newcommand{\proof}{{\it Proof.}\ } 
\newcommand{\qed}{\hfill  $\Box$\separation} 
 \def\seven{Z}
 \def\six{{\cal S}}

\begin{document}

\title{On the hypersurface of L\"uroth quartics}
 \author{Giorgio Ottaviani - Edoardo Sernesi\footnote{Both authors are members of GNSAGA-INDAM.}} 
\date{} 
 \maketitle

\abstract{The hypersurface of Luroth quartic curves inside the projective space of
plane quartics has degree 54. We give a proof of this fact along the lines
outlined in a paper by Morley, published in 1919. Another proof has been given
by Le Potier and Tikhomirov in 2001, in the setting of moduli spaces of vector
bundles on the projective plane. Morley's proof uses the description of plane
quartics as branch curves of Geiser involutions and gives new geometrical
interpretations of the 36 planes associated to the Cremona
hexahedral representations of a nonsingular cubic surface. }
 
\section*{Introduction}

In his  celebrated paper  \cite{jL68} L\"uroth proved that  a   nonsingular  quartic plane curve,  containing the ten vertices of a complete pentalateral, contains infinitely many such 10-tuples. This implies that such curves, called \emph{L\"uroth quartics}, fill an open set of an irreducible, SL$(3)$-invariant,  hypersurface $\L\subset\P^{14}$.  In his short paper \cite{fM14} Morley  computed  the degree 
of the L\"uroth hypersurface $\L$ by introducing some interesting ideas which seem to have been  forgotten, maybe because a few arguments are somehow obscure. 
 In this paper we put Morley's result and method on a solid foundation by  reconstructing his proof as faithfully as possible. The main result is the following:
 
 \begin{theorem}\label{T:main0}
 The L\"uroth hypersurface $\L \subset \P^{14}$ has degree 54.
 \end{theorem} 

Morley's proof uses the description of plane quartics as branch curves of  the degree-two rational self-maps of $\P^2$ called Geiser involutions. Every such involution is determined by the linear system of cubics having as base locus a 7-tuple of distinct points
$\seven=\{P_1, \dots, P_7\}$;  let's denote by $B(\seven) \subset \P^2$ the corresponding quartic branch curve. He introduces a closed condition on the space of such 7-tuples, given by the vanishing of the pfaffian of a natural skew-symmetric bilinear form between conics associated to each such $\seven$. By this procedure one obtains an irreducible  polynomial $\Psi(P_1,\dots, P_7)$  multihomogeneous of degree three in the coordinates of the points $P_1,\dots, P_7$, and skew-symmetric with respect to their permutations. We call $\Psi$   \emph{the Morley invariant}. The symbolic expression of $\Psi$ is  related to $\P^2_{\Z/2\Z}$, classically known as 
\emph{the Fano plane} (see \S \ref{S:cremona}).

\noindent
Then Morley proceeds in proving that the nonsingular quartics $B(\seven)$ corresponding to the 7-tuples $\seven$ for which the Morley invariant vanishes are precisely the L\"uroth quartics. This step of the proof uses a result of Bateman \cite{hB14} which gives an explicit description of an irreducible   13-dimensional family of configurations $\seven$ such that $B(\seven)$  is L\"uroth: Morley  shows that the Bateman configurations  are precisely those making $\Psi$ vanish.     
In order to gain control on the degree of $\L$  one must  consider    the full locus  of configurations $\seven$ such that $B(\seven)$ is a L\"uroth quartic, which  contains the locus  of      Bateman  configurations as a component. This can be realized as follows. Fix  six general points $P_1,\dots,P_6 \in \P^2$: the condition
$\Psi(P_1,\dots, P_6,P_7)=0$
on the   point $P_7$ defines  a plane cubic $E_{P_1,\dots, P_6}$ containing $P_1,\dots, P_6$, thus corresponding to  a plane section   $S \cap \Xi$ of the cubic surface $S \subset \P^3$ associated to the linear system of plane cubics through $P_1,\dots, P_6$.  
The plane $\Xi$ can be described explicitly  by means of the invariants introduced by Coble and  associated to the Cremona 
hexahedral equations of   $S$;  we   call  $\Xi$ the \emph{Cremona plane}. 
By construction the branch curve of the projection of $S$ to $\P^2$ from a point of $S \cap \Xi$ is L\"uroth. Conversely, given a general cubic surface $S \subset \P^3$ we obtain as many such plane sections as  the number of double-sixes on $S$, i.e. 36. 
   The final part, which relates the numbers 36 and 54, was implicitly  considered to be well-known by Morley.  We have supplied a proof which uses vector bundle techniques (see Theorem \ref{T:54}).

In order to put our work in perspective it is worth recalling here   some recent work related to L\"uroth quartics.  Let   $M(0,4)$ be the   moduli space of stable rank-two vector bundles on ${\bf P}^2$   with $(c_1,c_2)=(0,4)$.
Let 
\[
J(E)=\left\{l\in\left({\bf P}^2\right)^{\vee}|E_{l}\neq{\cal O}_l^2\right\}
\]
 be its curve of jumping lines. Barth proved in \cite{Bar} the remarkable facts that $J(E)$ is a L\"uroth quartic and that $\dim [M(0,4)]=13$. The Barth map, in this case, is the morphism
$$
\xymatrix{b: M(0,4) \ar[r]& \P^{14}, & E\ar@{|->}[r]&J(E)}
$$

It is well known that $b$ is generically finite and moreover that
\begin{equation}\label{E:barth}
\deg(b)\cdot \deg[\Im(b)] = 54 
\end{equation}
Indeed the value $54$ corresponds to the Donaldson invariant $q_{13}$ of ${\bf P}^2$ and it has been computed
by Li and Qin in \cite{LQ} theor. 6.29, and independently by Le Potier, Tikhomirov and Tyurin, see \cite{EPS} and the references therein. Another proof, related to secant varieties,  is in   Theorem 8.8 of \cite{gO07}.
Thanks to the result of Barth  mentioned above, the (closure of the) image of $b$ can be identified with the L\"uroth hypersurface $\L$, and   Theorem \ref{T:main0}, originally due to Morley, implies that
$\deg[\Im(b)] = 54 $.
The obvious corollary   is that
 $\deg (b)=1$, 
 that is \emph{the Barth map $b$ is generically injective}.
  This last result was obtained by Le Potier and Tikhomirov in \cite{PT} with a subtle and technical  degeneration argument. It also implies Theorem \ref{T:main0} via the identity (\ref{E:barth}).
 Our approach, which  closely follows  \cite{fM14}, is   more elementary and direct.  Le Potier and Tikhomirov also proved the injectivity of the Barth map for all higher values of $c_2$, treating  the case $c_2=4$ as the starting point of their inductive argument.
  
\separation

{\bf Acknowledgements.}  We thank A. Conca for kindly sharing his insight on syzygies of finite sets of points in $\P^2$. We are also grateful to I. Dolgachev for his encouragement  and for several remarks which helped us to improve the original version of this paper.

  \section{Apolarity}\label{S:apol}
 
We will work over an algebraically closed field $\k$ of characteristic zero.  
 Let $V$ be   a \k-vector space   of dimension 3 and denote by $V^\vee$ its dual.  
 The canonical bilinear form
 $\xymatrix{V \times V^\vee \ar[r] & \k}$
 extends to a natural pairing:
 \begin{equation}\label{E:pol1}
 \xymatrix{S^dV \times S^nV^\vee \ar[r] & S^{n-d} V^\vee}, 
 \qquad \xymatrix{(\Phi,F)\ar@{|->}[r]&P_\Phi(F)}
 \end{equation}
 for each $n \ge  d$, which is
   called \emph{apolarity}.  $\Phi$ and $F$ will be called \emph{apolar} if  $P_\Phi(F)=0$.
 
\noindent
    After choosing a basis of $V$ we can identify  the symmetric algebra $\mathrm{Sym}(V^\vee)$ with the polynomial algebra
    $\k[X_0,X_1,X_2]$ and $\mathrm{Sym}(V)$ with $\k[\partial_0,\partial_1,\partial_2]$, where 
 $\partial_i :={\partial\over \partial X_i}$, $i=0,1,2$, are the \emph{dual indeterminates}. With this notation   apolarity is the natural pairing between differential operators and polynomials. We can also identify 
 $\P(V)=\P^2$ and $\P(V^\vee)=\P^{2\vee}$.   
 
 \noindent
  Elements of $S^dV$, up to a non-zero factor, are called \emph{line curves} of degree $d$ (line conics, line cubics, etc.) while elements of $S^dV^\vee$, up to a non-zero factor,  are \emph{point curves} of degree $d$ (point conics, point cubics, etc.). 
  We will be mostly interested in the case of degree $d=2$.  In this case, in coordinates, apolarity  takes the form:
\[
 P_\Phi(\sum_{ij}\alpha_{ij}X_iX_j)=  \sum_{ij} a_{ij}\alpha_{ij}
\]
if $\Phi = \sum_{ij}a_{ij}\partial_i\partial_j$.
 Suppose given a  point conic defined by the polynomial
\begin{equation}\label{E:apol0}
\theta=  \sum_{ij} A_{ij}X_iX_j \in S^2V^\vee
\end{equation}
Assume that $\theta$ is \emph{nonsingular}, i.e. that its coefficient  matrix  $(A_{ij})$ is invertible, and     consider its \emph{dual curve}    $\theta^*=\sum_{ij}a_{ij}\partial_i\partial_j$.
\noindent
We will say that \emph{a point conic 
$C= \sum_{ij} \alpha_{ij}X_iX_j \in S^2V^\vee$
is conjugate to} $\theta$ if it is apolar to $\theta^*$. This gives a notion of conjugation between point conics and, dually, between line conics. Note that if $C$ is conjugate to $\theta$ then it is not necessarily true that $\theta$ is conjugate to $C$, i.e. this notion is not symmetric. In particular,  we did not require $C$ to be nonsingular in the definition.  

Another important special case of (\ref{E:pol1}) is the following. Given a point $\xi=(\xi_0,\xi_1,\xi_2) \in \P^2$, the corresponding linear form 
$\xi_0\partial_0+\xi_1\partial_1+\xi_2\partial_2 \in (V^\vee)^\vee=V$ will be also denoted by $\Delta_\xi$ and called the \emph{polarization operator} with pole $\xi$. For each $d \ge 2$   it defines a linear map:
\[
\begin{array}{ll}
\xymatrix{\Delta_\xi:S^dV^\vee \ar[r]& S^{d-1}V^\vee} \\
\Delta_\xi F(X) = \xi_0\partial_0F(X) +\xi_1\partial_1F(X) 
+\xi_2\partial_2F(X) 
\end{array}
\]
associating to a homogeneous polynomial $F$ of degree $d$  a homogeneous  polynomial of degree $d-1$ called \emph{the (first) polar} of $\xi$ with respect to $F$.  
Higher polars are defined similarly by iteration.  

\noindent
Consider  the case $d=2$. 
Given a nonsingular point conic $\theta$, polarity associates to each point $\xi\in \P^2$, the \emph{pole},  its polar line  $\Delta_\xi\theta$, and this gives an isomorphism $\P^2 \cong \P^{2\vee}$. 
  Two points   will be called \emph{conjugate   with respect to $\theta$} if each of them belongs to the polar line of  the other.   Two lines are called \emph{conjugate with respect to} $\theta$ if each of them contains the pole of the other. We will need the following elementary properties of apolarity, whose proof we leave to the reader.

  \begin{prop}\label{P:apol1}  Let $\theta$ be a nonsingular point conic.
  \begin{itemize}
  \item[(i)]  A point conic $C$ reducible in two distinct lines $\ell_1\ell_2$ is conjugate to  $\theta$ if and only if the two lines are conjugate with respect to $\theta$ or, equivalently, if and only if $\ell_1$ and $\ell_2$  are conjugate in the involution on the pencil of lines through the point $\ell_1 \cap \ell_2$ having as fixed    points the tangent lines to $\theta$.
  
  \item[(ii)]   A point conic  $C$ consisting of a double line  is conjugate to  $\theta$ if and only if the line is tangent to $\theta$.

  \item[(iii)]  Every point conic $C$ reducible in the tangent line to $\theta$ at a point 
  $\xi \in \theta$ and in any other line through $\xi$ is conjugate to $\theta$.

  \item[(iv)]  If a point conic C is conjugate to $\theta$ then for each $\xi \in C$ the reducible point conic consisting of the lines joining $\xi$ with $C \cap \Delta_\xi\theta$ is also conjugate to $\theta$.
  
  \end{itemize}
  
   \end{prop}

 Given a nonsingular point conic $\theta$, we will call a point cubic $D \in S^3V^\vee$  \emph{apolar to}  $\theta$ if  $P_{\theta^*}(D) = 0$, i.e. if $\theta^*$ and $D$ are apolar.  Note that, since $P_{\theta^*}(D) \in V^\vee$, the condition of apolarity to $\theta$ is equivalent to  three linear conditions on point cubics. We will need the following.
 
 \begin{prop}\label{P:apol2}
 Let $\theta$ be a nonsingular point conic and $D$ a point cubic apolar to $\theta$. Then for every point $\xi \in \P^2$ the polar conic $\Delta_\xi D$ is conjugate to $\theta$.
 \end{prop}
 
 \proof
 From $P_{\theta^*}(D) = 0$ it follows that for any $\xi\in \P^2$ we have:
 \[
 0 = P_{\Delta_\xi\theta^*}(D) = P_{\theta^*\Delta_\xi}(D) =  P_{\theta^*}(\Delta_\xi D)
 \]
 \qed

 \begin{prop}\label{P:apol3}
 Let $\theta$ be a nonsingular point conic. Then:
 \begin{itemize}
 
 \item[(i)] For every line $L$ the point cubic $\theta L$ is not apolar to $\theta$.
 
 \item[(ii)]  For every effective divisor $\sum_{i=1}^6 P_i$ of degree six on $\theta$ there is a unique point cubic $D$ such that $D$ is apolar to $\theta$ and $D \cdot \theta = \sum_{i=1}^6 P_i$.  If  $\sum_{i=1}^6 P_i$ is general then $D$ is irreducible. 
 
 \item[(iii)]  For every effective divisor $\sum_{i=1}^5 P_i$ of degree five on $\theta$ and for a general point 
 $P_6 \notin \theta$ there is a unique point cubic $D$ containing $P_6$ such that $D$ is apolar to $\theta$ and $D \cdot \theta > \sum_{i=1}^5 P_i$.

 \end{itemize}
 \end{prop}
 
 \proof   
 (i)   Let  $\xi\in \theta$ but $\xi\notin L$. Then:
  \[
  \begin{array}{l}
 \Delta_\xi[P_{\theta^*}(\theta L)] = P_{\theta^*}(\Delta_\xi(\theta L))=
  P_{\theta^*}[\Delta_\xi(\theta)L+\theta \Delta_\xi(L)]   \\
  = P_{\theta^*}(\Delta_\xi(\theta)L) + P_{\theta^*}(\theta \Delta_\xi(L)) = 0+3\Delta_\xi(L) \ne 0
 \end{array}
  \]
 where the last equality is by Prop. \ref{P:apol1}(iii). Therefore $P_{\theta^*}(\theta L)\ne 0$.

 (ii) If  $C$ is a point cubic such that $C\cdot \theta = \sum_{i=1}^6 P_i$ then all other point cubics with this property are of the form $D=C- \theta L$ for some line $L$.  Taking $\xi\in \theta$, by the previous computation we obtain:
 \[
 \Delta_\xi(P_{\theta^*}D)=\Delta_\xi\left[ P_{\theta^*}(C-\theta L)\right] = \Delta_\xi P_{\theta^*}(C) - \Delta_\xi(3L)
 \]
 This is zero for all $\xi\in\theta$ if and only if  $P_{\theta^*}(D)=0$, if and only if
 $3L=P_{\theta^*}(C)$.
 Finally, the six-dimensional linear system of cubics apolar to $\theta$ cannot consist  of reducible cubics.

 (iii) follows easily from (ii). \qed
    
We refer the reader to \cite{D} for a more detailed treatment  of polarity and apolarity. From now on by a \emph{conic,}  resp a \emph{cubic}, etc.,   we will mean a point conic (resp. point cubic, etc.) unless otherwise specified.

  \section{The Morley form}\label{S:morley}

 Consider   seven distinct points $P_1, \dots, P_7 \in \P^2$ and let $\seven =\{P_1, \dots, P_7\}$.  Let   $\I_\seven \subset  {\cal O}_{\P^2}$ be the ideal sheaf of $\seven$ and 
 \[
 I_\seven = \bigoplus_k I_{\seven,k}=\bigoplus_kH^0(\P^2, \I_\seven(k)) \subset \k[X_0,X_1,X_2]
 \]
  the homogeneous ideal of $\seven$.

 \begin{prop}\label{L:alg1}
 Assume that $\seven$ is not contained in a conic.  Then: 
\begin{itemize}
\item[(i)] There is a matrix of homogeneous polynomials
\[
 A= \pmatrix{L_0(X)&L_1(X)&L_2(X) \cr \theta_0(X)&\theta_1(X) & \theta_2(X)}
  \]
  where  $\deg(L_i(X))=1$ and $\deg(\theta_i(X))=2$ for $i=0,1,2$ such that  $I_\seven$ is generated by the maximal minors of $A$.  

\item[(ii)] Six of the seven points  $P_1, \dots, P_7$ are on a conic  if and only if for any matrix $A$ as in (i) the linear forms  $L_0(X),L_1(X),L_2(X)$ are linearly dependent. 

\end{itemize}
  \end{prop}

 \proof (i) By the Hilbert-Burch theorem the homogeneous ideal of any finite set of points in $\P^2$ is generated by the maximal minors of a $t\times (t+1)$ matrix $A$ of homogeneous polynomials of positive degrees for some $t \ge 1$ (\cite{dE05}, Th. 3.2). 
Since $\seven$ is contained in at least three linearly independent cubics it must be $t \le 2$. Since moreover $\seven$ is not a complete intersection of two curves, we must have $t=2$.  The numerical criterion of \cite{CGO}  (see also  \cite{dE05}, Cor. 3.10) shows that the only possibility is the one stated. 

\noindent
(ii) Clearly it suffices to prove the assertion for one matrix $A$ as in (i). Assume that $P_1,\dots,P_6$ are on a conic $\theta_0$, and that $P_7 \notin \theta_0$. Let $L_1,L_2$ be two distinct lines through $P_7$. Then $\langle C_0, L_2\theta_0,-L_1\theta_0\rangle= H^0(\I_\seven(3))$ for some cubic $C_0$, and since $P_7 \in C_0$ there are conics
$\theta_1,\theta_2$ such that
 $C_0=L_1\theta_2 - L_2\theta_1$
 so that  we can take
\begin{equation}\label{E:A}
 A= \pmatrix{0&L_1(X)&L_2(X) \cr \theta_0(X)&\theta_1(X) & \theta_2(X)}
  \end{equation}
 and $L_0=0,L_1(X),L_2(X)$ are linearly dependent. 
 Conversely, assume that $L_0(X),L_1(X),L_2(X)$ are linearly dependent for some $A$ as in (i).  After multiplying to the right by a suitable element of SL$(3)$  we may assume that $L_0(X)=0$, i.e. that $A$ has the form (\ref{E:A}). It immediately follows that one of the seven points is  $L_1 \cap L_2$ and that the other six are contained in $\theta_0$.   \qed

  Unless otherwise specified, \emph{we will always assume that 
  $\seven = \{P_1,\dots,P_7\}$  consists of distinct points  not on a conic.}

\noindent
  Let  $\xi = (\xi_0,\xi_1,\xi_2)$ be new indeterminates, and consider the polynomial:
  \[
\begin{array}{ll}
  S(\xi,X) := \left|\matrix{L_0(\xi)& L_1(\xi)&L_2(\xi) \cr 
L_0(X)&L_1(X)&L_2(X) \cr \theta_0(X)&\theta_1(X) & \theta_2(X)} 
\right|   \matrix{=L_0(\xi)C_0(X)+L_1(\xi)C_1(X)\cr +L_2(\xi)C_2(X)}
 \end{array}
 \]
  where the $L_j$'s and the $\theta_j$'s are the entries of a matrix $A$ as in (i) of Proposition \ref{L:alg1}.  $S(\xi,X)$ is bihomogeneous of degrees 1 and 3 in $\xi$ and $X$ respectively.   
  
  \noindent
Given points $P=(x_0,x_1,x_2), Q=(y_0,y_1,y_2) \in \P^2$ we will denote by
  \[
  |PQX| = \left|\matrix{x_0&x_1&x_2 \cr y_0&y_1&y_2 \cr X_0&X_1&X_2 } \right|
  \]
   If $P \ne Q$ then  $ |PQX| =0$ is the line containing $P$ and $Q$. 
     
  \begin{lemma}\label{L:alg2}
  Let $\seven$ and $S(\xi,X)$ be as above. Then:
\begin{itemize}
  \item[(i)] Up to a constant factor $S(\xi,X)$ depends only on $\seven$ and not on the particular  choice of the matrix $A$. 
  \item[(ii)] If no six of the points of $\seven$ are on a conic then for every choice of $\xi \in \P^2$ the cubic  $S(\xi,X)$ is not identically zero, contains $\xi$ and $\seven$ and  is singular at $\xi$ if   $\xi \in \seven$.  All cubics in $H^0(\P^2,\I_\seven(3))$ are obtained as $\xi$ varies in $\P^2$.
  
  \item[(iii)]  If $\{P_1,\dots, P_6\}$ are on a nonsingular conic $\theta$ and $P_7 \notin \theta$ then  
    \[
    S(\xi,X) = |P_7\xi X| \theta
    \]
  In particular $S(P_7,X) \equiv 0$
    and only the 2-dimensional vector space of reducible cubics in $H^0(\I_\seven(3))$ is represented in the form $S(\xi,X)$. 
    
        \end{itemize}
    \end{lemma}
    
    \proof
    (i) A different choice of the matrix $A$ can    be obtained by multiplying it  on the right by some $M \in {\rm GL}(3)$ and this has the effect of changing $S(\xi,X)$ into $S(\xi,X){\rm det}(M)$. Also left action is possible but it does not change $S(\xi,X)$.

 \noindent
 (ii) Since $L_0(\xi),L_1(\xi),L_2(\xi)$ are linearly independent (Lemma \ref{L:alg1})   $S(\xi,X)$ cannot be identically zero, and it follows that all of $H^0(\P^2,\I_\seven(3))$ is obtained in this way. Clearly $S(\xi,X)$ contains $\xi$. From the identity:
 \[ 
 0 =  {\partial \left[\sum_j L_j(X)C_j(X)\right] \over \partial X_h} = 
 \sum_j{\partial L_j(X)\over \partial X_h}C_j(X) + \sum_j L_j(X) {\partial C_j(X) \over \partial X_h}
 \]
 we deduce: 
 \[
 {\partial S(\xi,X) \over \partial X_h} = \sum_jL_j(\xi){\partial C_j(X) \over \partial X_h} =
 -\sum_j{\partial L_j(\xi)\over \partial \xi_h}C_j(X)
 \]
 The last expression for the partials of  $S(\xi,X)$ shows that  
 \[
  {\partial S(\xi,X) \over \partial X_h}(\xi) = 0
  \]
  for $h=0,1,2$  if $\xi\in\seven$,  so that 
 $S(\xi,X)=0$ is singular at $\xi$ in this case. 
 
 \noindent
 (iii) As in the proof of Lemma \ref{L:alg1}, we can choose $L_0=0$, and $L_1$ and $L_2$ linearly independent and containing $P_7$ and $\{C_0, L_2(X)\theta,-L_1(X)\theta \}$ as a basis of 
 $H^0(\P^2,\I_\seven(3))$. Then
  \[
 S(\xi,X) = [L_1(\xi)L_2(X)-L_2(\xi)L_1(X)]\theta
 \]
 From this expression (iii) follows immediately.
 \qed

 Lemma \ref{L:alg2} shows that $S(\xi,X)$ is uniquely determined by $\seven$ up to a constant factor. More precisely we have the following:
 
 \begin{prop}\label{P:deg}
 The coefficients of $S(\xi,X)$ can be expressed as   multihomogeneous  polynomials  of degree 5  in the   coordinates of the points $P_1,\dots, P_7$ which are symmetric with respect to permutations of $P_1,\dots, P_7$. 
 \end{prop}

 \proof
 On  $\P^2\times \P^2$ with homogeneous coordinates $\xi$ and $X$ consider 
 the exact sequence:
 \[
   0 \to \I_\Delta(1,3) \to  {\cal O}_{\P^2\times\P^2}(1,3) \to  {\cal O}_\Delta(4) \to 0
    \]
  where $\Delta \subset\P^2\times \P^2$  is the diagonal.  Since
  $h^0( {\cal O}_{\P^2\times\P^2}(1,3))=30$ and $h^0( {\cal O}_\Delta(4))=15$
    from the exact sequence we deduce that
 $h^0(\I_\Delta(1,3))=15$
 and that  $S(\xi,X) \in H^0(\I_\Delta(1,3))$.  Given a polynomial:
 \[
 P(\xi,X) = \sum_j \xi_j D_j(X)  \in H^0( {\cal O}_{\P^2\times\P^2}(1,3))
 \]
 the condition that it belongs to  $H^0(\I_\Delta(1,3))$ corresponds to the vanishing of the 15 coefficients of $P(X,X) \in H^0( {\cal O}_\Delta(4))$, and these are 15  linear homogeneous conditions with constant coefficients on the 30 coefficients of $P(\xi,X)$. The condition that $P(\xi,X)=S(\xi,X)$ up to a constant factor is that   moreover:
 \begin{equation}\label{E:deg1}
 \sum_j\xi_jD_j(P_i) =0, \qquad  i=1,\dots, 7
 \end{equation}
 because this means that  the cubic $P(\xi,X)=0$ contain $\seven$ for all $\xi\in \P^2$. For each $i=1,\dots,7$ the condition (\ref{E:deg1}) means
 \begin{equation}\label{E:deg2}
 D_0(P_i) = D_1(P_i) = D_2(P_i) = 0
 \end{equation}
 and these are 3 linear homogeneous conditions on the 30 coefficients of $P(\xi,X)$ with coefficients which are homogeneous of degree 3 in $P_i$. Since
 $P(\xi,X) \in H^0(\I_\Delta(1,3))$, we also have:
 \begin{equation}\label{E:deg3}
 \sum_j x_{ij} D_j(P_i) = 0
 \end{equation}
 where $P_i = (x_{i0},x_{i1}, x_{i2})$.  This condition implies that  only 2 of the 3 conditions (\ref{E:deg2}) are independent: if say $x_{i0}\ne 0$ then we can choose $D_1(P_i) = D_2(P_i) = 0$. Moreover whenever either one is satisfied,   the remaining one is divisible by $x_{i0}$,  thanks to the relation (\ref{E:deg3}).  Therefore for each $i=1,\dots, 7$ we obtain two linear homogeneous conditions on the 30 coefficients of $P(\xi,X)$, with coefficients which are homogeneous of degree 3 and 2 respectively in $P_i$. 
  Altogether we obtain $29 = 15+14$ linear homogeneous conditions on the 30 coefficients of $P(\xi,X)$.  The maximal minors of  their coefficient matrix are the coefficients of $S(\xi,X)$, and they are  multihomogeneous of degree 5 in the $P_i$'s by what we have shown.  Any transposition of 
 $P_1,\dots, P_7$ permutes two pairs of adjacent  rows of the matrix so that the maximal minors remain unchanged. \qed

  \begin{definition}\label{D:mor}
 The \emph{Morley form} of $\seven$   is the biquadratic homogeneous polynomial
   in $\xi,X$:
 \[
 M(\xi,X) := \Delta_\xi S(\xi,X) = 
\left|\matrix{L_0(\xi)& L_1(\xi)&L_2(\xi) \cr 
L_0(X)&L_1(X)&L_2(X) \cr \Delta_\xi \theta_0(X)&\Delta_\xi \theta_1(X) &\Delta_\xi  \theta_2(X)} 
\right|
 \]
 \end{definition}
 
  For every $\xi \in \P^2$ such that $S(\xi,X)$ is not identically zero $M(\xi,X)$  represents in $\P^2=\Proj(\k[X_0,X_1,X_2])$ the polar conic of $\xi$ with respect to the cubic  $S(\xi,X)$. 
  Clearly
 it contains $\xi$ and, if $\xi \in \seven$, it is reducible into the principal tangent lines at $\xi$ of 
    $S(\xi,X)$ by Lemma \ref{L:alg2}.  Note that, by Lemma \ref{L:alg2}, $S(\xi,X)\equiv 0$ (and consequently 
    $M(\xi,X)\equiv 0$)  if and only if six of the seven points of $\seven$ are on a conic and $\xi$ is the seventh point. 
 Since $M(\xi,\xi)=0$,  the Morley form is skew-symmetric in $\xi,X$. Therefore its $6\times 6$ matrix of coefficients  $(M_{hk})$ has  determinant which  is the square of its pfaffian.

 \begin{corollary}\label{C:deg1}
 The pfaffian of $M(\xi,X)$ can be expressed as a polynomial $F(P_1, \dots, P_7)$  multihomogeneous  of degree 15 in the coordinates of the points 
 $P_1,\dots, P_7$ and symmetric with respect to permutations of the points. 
  \end{corollary}

  \proof By Proposition \ref{P:deg} the coefficients   $M_{hk}$ are  multihomogeneous of degree 5 in the coordinates of the $P_i$'s.  Therefore the determinant is  multihomogeneous of degree 
  $30$ in the $P_i$'s, thus the pfaffian has degree 15 in each of them.  The symmetry follows from that of the coefficients  $M_{hk}$ which holds by Proposition \ref{P:deg} . \qed

 The Morley form $M(\xi,X)$ defines a bilinear skew-symmetric form:
 \[
 \xymatrix{S^2V^\vee \times S^2V^\vee \ar[r] & \k}
 \]
  If $F(P_1, \dots, P_7)=0$ then this form is degenerate.  The 7-tuples  $\{P_1,\dots, P_7\}$ of points in $\P^2$ for which this happens are such that, when $\xi$ varies in $\P^2$,     all the conics $M(\xi,X)$ are contained in a hyperplane of 
  $\P(S^2V^\vee)$.  The search for such 7-tuples is our next goal.

    \section{The Morley invariant} \label{S:morleinv}

    \begin{prop}\label{P:apolar1}
    If  $\seven=\{P_1,\dots, P_7\}$ consists of distinct points not on a conic, six of which are on a conic, then   $F(P_1, \dots, P_7)=0$. 
     \end{prop}
    
    \proof 
   Let $\theta$ be the conic containing six of the seven  points, say $P_1,\dots, P_6$.  From Lemma \ref{L:alg2}(iii)   it follows that 
$S(\xi,X) = \theta |P_7\xi X|$. 
Therefore  all the conics $M(\xi,X)$, $P_7\ne \xi\in \P^2$,     are   contained in the hyperplane   $H_{P_7}\subset S^2V^\vee$ of conics which contain $P_7$.  This implies that the skew-symmetric form
   $M: S^2V^\vee \times S^2V^\vee \to \k$
   is degenerate, hence its pfaffian vanishes. \qed

   Given $p_1,\dots,p_6 \in  \P^2$,    define as in \cite{Cob2} pag. 136 (see also \cite{DO88}, p. 191):
 \[
{\cal Q}(p_1,\dots,p_6) = |134| |156|  |235| |246|- |135| |146| |234|  |256|
\] 
where we   use the symbolic notation:
\[
|ijk| := \left|\matrix{p_{i0}&p_{i1}&p_{i2} \cr p_{j0}&p_{j1}&p_{j2} \cr p_{k0}&p_{k1}&p_{k2}}\right|
\]
${\cal Q}(p_1,\dots,p_6)$ is a  multihomogeneous polynomial of degree two
in the coordinates of the points $p_1,\dots,p_6$,  skew-symmetric w.r. to
them and which vanishes if and only if $p_1,\dots,p_6$ are on a conic.
Moreover ${\cal Q}(p_1,\dots,p_6)$ is irreducible because any factorization
would involve invariants of lower degree for the group ${\rm SL}(3)\times
{\rm Alt}_6$ which do not exist.

\begin{prop}\label{P:pfaff2}
Consider  distinct  points $P_1, \dots, P_7$.  The polynomial 
\begin{equation}\label{E:prod1}
\prod_i {\cal Q}(P_1, \dots, \widehat P_i, \dots, P_7)  
\end{equation}
 is  multihomogeneous of degree 12 in the coordinates of each point $P_i$,  $i=1,\dots,7$ and  skew-symmetric w.r. to permutations of   $P_1, \dots, P_7$. 
 It vanishes precisely on the 7-tuples which contain six points on a conic    and divides the pfaffian polynomial $F(P_1,\dots,P_7)$.
\end{prop}

\proof 
 Since each polynomial
${\cal Q}(P_1, \dots, \widehat P_i, \dots, P_7)$
has degree two in the coordinates of each of the six points $P_1, \dots, \widehat P_i, \dots, P_7$ it follows that  the product (\ref{E:prod1})   is  multihomogeneous of degree 12 in the coordinates of each point $P_i$,  $i=1,\dots,7$.
Let $1 \le   i<j \le 7$.  Then  each  ${\cal Q}(P_1, \dots, \widehat P_k, \dots, P_7)$,  $k \ne i,j$, is skew-symmetric w.r. to $P_i$ and $P_j$. 
On the other hand    
\[
{\cal Q}(P_1, \dots, \widehat P_i, \dots, P_7){\cal Q}(P_1, \dots, \widehat P_j, \dots, P_7)
\]
 is symmetric w.r. to $P_i$ and $P_j$ because
\[
{\cal Q}(P_1,\dots,\widehat P_i,\dots,P_i,\dots,P_7) = (-1)^{i-j+1} {\cal Q}(P_1, \dots, \widehat P_j, \dots, P_7)
\]
where on the left side $P_j$ has been replaced by $P_i$ at the $j$-th place.  Therefore (\ref{E:prod1}) is skew-symmetric.
 It is clear that (\ref{E:prod1}) vanishes precisely at those 7-tuples which include six points on a conic.
 The last  assertion  follows at once from   Proposition \ref{P:apolar1}.  \qed

    We will denote by $\Psi(P_1,\dots,P_7)$ the polynomial such that 
    \begin{equation}\label{E:pfaff1}
    F(P_1,\dots,P_7) =   \Psi(P_1,\dots,P_7) \prod_i {\cal Q}(P_1, \dots, \widehat P_i, \dots, P_7) 
    \end{equation}
    We call $\Psi(P_1,\dots,P_7)$  the \emph{Morley invariant} of the seven points  $P_1,\dots,P_7$.  
    
    \begin{corollary}\label{C:pfaff1}
    The Morley invariant $\Psi(P_1,\dots,P_7)$ is  multihomogeneous of degree 3 in the coordinates of the points $P_i$ and skew-symmetric with respect to  $P_1,\dots,P_7$.
    \end{corollary}
    
    \proof
    It follows   from Corollary \ref{C:deg1} and from the fact that  the polynomial (\ref{E:prod1})  is  multihomogeneous of degree 12 and skew-symmetric. \qed

  Let $\seven = \{P_1,\dots, P_7\}$ be given  consisting of distinct points   not on a conic as always.  The  net of cubic curves $|H^0(\I_\seven(3))|$  contains a unique cubic singular at $P_i$  for each $i=1,\dots, 7$:  we denote by  $M^\seven_{P_i}$, or simply by 
          $M_{P_i}$ when no confusion is possible, the reducible conic  of its principal tangents at $P_i$.  
          
          \noindent
          As already remarked after Definition \ref{D:mor}, 
           if no six of the points of $\seven$ are on a conic  then $M_{P_i}=M(P_i,X)$ for all $1 \le i \le 7$.     
  If instead     six of the points, say  $P_1,\dots, P_6$,  are on a conic $\theta$  
   then  $M_{P_i}=M(P_i,X)$ for   $i=1,\dots,6$, but $M_{P_7}$ \emph{ is not obtained from} $M(\xi,X)$. 
  By Lemma \ref{L:alg2} we have  
  $S(\xi,X) =  \theta |P_7\xi X|$
  and  therefore if $\xi \ne P_7$ then  $M(\xi,X)$   is reducible in the line $|P_7\xi X|$ and in the polar line of $\xi$ with respect to $\theta$.  
 
 \begin{prop}\label{P:apolar0}
 Assume that  $\seven = \{P_1,\dots, P_7\}$ are such that six of them, say  $P_1,\dots, P_6$,  are on a nonsingular conic $\theta$.  Then:
 
 \begin{itemize}
 
  \item[(i)]   the conics $M(\xi,X)$, as $\xi$ varies in $\P^2\backslash \{P_7\}$,  depend only on $\theta$ and $P_7$, and not on the points $P_1, \dots, P_6$.  They  generate a vector subspace of dimension 4 of $S^2V^\vee$ which is the intersection of the hyperplane $H_\theta$ of conics conjugate to $\theta$ with the hyperplane $H_{P_7}$ of conics containing $P_7$.  Moreover
  \begin{equation}\label{E:morleinv1}
  H_\theta \cap H_{P_7}  = \langle M_{P_1}, \dots, M_{P_6} \rangle
  \end{equation}
  for a general choice of $P_1, \dots, P_6 \in \theta$. 
  
   \item[(ii)]  for a general choice of  $P_1,\dots, P_6 \in \theta$ and of  $P_7 \notin \theta$  the reducible conic $M_{P_7}$  is not conjugate to $\theta$. In particular 
   \[
   \langle M_{P_1}, \dots, M_{P_6}, M_{P_7}\rangle = H_{P_7}
   \]
   has dimension 5.
      \end{itemize}
  \end{prop}

 \proof
  (i) We can assume that:
  $\theta=X_0^2+ 2 X_1X_2$
  and that $P_7 = (1,0,0)$. Then  $\theta X_1, \theta X_2 \in H^0(\I_\seven(3))$, and therefore:
  \[
  S(\xi,X) = \theta(\xi_1 X_2 - \xi_2X_1) 
  \]
  so that 
  \[
   M(\xi,X) = (\Delta_\xi \theta) (\xi_1 X_2 - \xi_2X_1) = 
   2(\xi_0X_0+\xi_1X_2+\xi_2X_1)(\xi_1 X_2 - \xi_2X_1)
  \]
Clearly this  expression does not depend on the points $P_1, \dots, P_6$. 

\noindent
The dual   of $\theta$ is the line conic
 $\theta^*=\partial_0^2+2\partial_1\partial_2$. 
  Therefore the hyperplane $H_\theta$ of conics conjugate to $\theta$ consists of the   conics 
  $C: \sum_{ij}\alpha_{ij}X_iX_j$ such that $\alpha_{00} + \alpha_{12} = 0$. The hyperplane $H_{P_7}$ of conics containing $P_7$ is given by the condition $\alpha_{00}=0$. Therefore $H_\theta \cap H_{P_7}$ has equations  $\alpha_{00} = \alpha_{12} = 0$. 
   Since $M(\xi,X)$ does not contain the terms $X_0^2$ and  $X_1X_2$, it follows that  
  $M(\xi,X) \in H_\theta \cap H_{P_7}$ for all $\xi \ne P_7$. 
  Now observe that:
  \[
  M(\xi,X) = \cases{X_1^2& if $\xi=(0,0,1)$ \cr 
  X_2^2 & if $\xi = (0,1,0)$ \cr 
  X_0X_1-X_1^2& if $\xi = (1,0,1)$ \cr
  X_0X_2+X_2^2& if $\xi = (1,1,0)$}
  \]
  which are linearly independent: it follows that the conics  $M(\xi,X)$  generate  $H_\theta \cap H_{P_7}$.  
  
  \noindent
  (\ref{E:morleinv1}) can be proved by a direct computation as follows. The conics $M(\xi,X)$ corresponding to the points 
  \[
  \xi = (0,0,1), (0,1,0), (2,-2,1),  (2i,-2i,-1) \in \theta
  \]
   are respectively:
  \[
  X_1^2, X_2^2, 2X_0X_1+4X_0X_2+X_1^2-4X_2^2, -2iX_0X_1+4iX_0X_2+X_1^2-4X_2^2
  \]
  and they are linearly independent.

  (ii) Keeping the same notations as above,  observe that a reducible conic   with double point $P_7$ and not conjugate to $\theta$ is of the form 
  $\alpha_{11}X_1^2 + \alpha_{22}X_2^2 + \alpha_{12}X_1X_2$
  for   coefficients  $\alpha_{11}, \alpha_{22}, \alpha_{12}$ such that  $\alpha_{12}\ne 0$. It follows that any cubic  of the form
 $D=X_0X_1X_2 + F(X_1,X_2)$
  where $F(X_1,X_2)$ is a general cubic polynomial,   is singular at $P_7$ and has the conic of principal tangents  equal to $X_1X_2$, and therefore not conjugate to $\theta$. Now it suffices to take 
  $\{P_1,\dots, P_6\} =  D \cap \theta$
  to have a configuration $\seven = \{P_1,\dots, P_7\}$ satisfying the desired conditions.  \qed
    
    \begin{corollary}\label{P:nonzero}
    The Morley invariant $\Psi(P_1,\dots,P_7)$ is not identically zero.
  \end{corollary}

    \proof  Since $M(\xi,X)$ is skew-symmetric,  for a given $\seven = \{P_1, \dots, P_7\}$ the subspace  $\Sigma_\seven \subset S^2V^\vee$ generated by the conics
    $M(\xi,X)$ when  $\xi$ varies in $\P^2$, has even dimension. Moreover, if  no six of the points of
      $\seven$ are on a conic  then
    $\langle M_{P_1}, \dots, M_{P_7}\rangle \subset \Sigma_\seven$.
   If moreover $P_1,\dots,P_7$ are general points then the space on the left hand side has dimension  $\ge 5$   because this happens for the special choice of $P_1, \dots, P_6$ on a nonsingular conic  and $P_7$ general (Prop. \ref{P:apolar0}).   Therefore we conclude that 
   $\Sigma_\seven = S^2V^\vee$ if $P_1,\dots,P_7$ are general points, and this means that the skew-symmetric form $M(\xi,X)$ is non-degenerate, equivalently its pfaffian does not vanish, and, a fortiori, 
   $\Psi(P_1,\dots,P_7) \ne 0$.  \qed

   \begin{remark}\label{R:mor1}\rm
  From Proposition \ref{P:nonzero} it follows that $\Psi$ defines a     hypersurface 
     \[
          V(\Psi) = \left\{(P_1,\dots,P_7) \in (\P^2)^7:  \Psi(P_1, \dots, P_7) = 0 \right\}
         \]
   in the 7-th cartesian product of $\P^2$. Since $\Psi$ is not divisible by 
   \[
   \prod_i {\cal Q}(P_1, \dots, \widehat P_i, \dots, P_7)
   \]
    the general element $(P_1,\dots,P_7)$  of each irreducible component of $V(\Psi)$ consists of points no six of which are on a conic. Moreover it follows from the proof of \ref{P:nonzero} that if 
    $(P_1,\dots,P_7)\in  V(\Psi)$ has such a property  then 
    $\dim \langle M_{P_1}, \dots, M_{P_7} \rangle  \le 4$.
 \end{remark}
 
Since $\Psi$ is skew-symmetric with respect to the action of $S_7$ on  $(\P^2)^7$, the hypersurface $V(\Psi)$ is $S_7$-invariant and therefore defines a hypersurface in the symmetric product   $(\P^2)^{(7)}$.

 \begin{definition}\label{D:W}
We will denote by ${\cal W}$ the image of  $V(\Psi)$ in the symmetric product
   $(\P^2)^{(7)}$.
\end{definition}

    \section{Cremona hexahedral equations}\label{S:cremona}
    
   In this section we want to show the relations between some classical work by  Cremona and Coble and  the  objects we have considered, and  to clarify their invariant-theoretic significance.  
   
   \noindent
   In  $\P^5$ with coordinates $(Z_0, \dots,Z_5)$  consider the following equations:
    \begin{equation}\label{E:cre1}
    \cases{Z_0^3+Z_1^3+Z_2^3+Z_3^3+Z_4^3+Z_5^3 = 0 \cr
    Z_0+Z_1+Z_2+Z_3+Z_4+Z_5 = 0  \cr
    \beta_0 Z_0+\beta_1 Z_1 + \beta_2 Z_2 + \beta_3 Z_3 + \beta_4 Z_4 + \beta_5 Z_5 = 0 }
    \end{equation}
    where the $\beta_s$'s are constants. These equations define a cubic  surface $S$ in a $\P^3$ contained 
    in $\P^5$. If  $S$ is nonsingular then (\ref{E:cre1}) are called \emph{Cremona hexahedral equations} of $S$, after \cite{Cre}.  They have several remarkable properties, the most important for us being that these equations also determine a double-six of lines on the surface $S$.  
    
    \noindent
    Recall that a \emph{double-six} of lines on a nonsingular cubic surface $S \subset \P^3$ consists of two sets of six lines
    $\Delta = (A_1, \dots, A_6; B_1, \dots B_6)$
    such that the lines $A_j$ are mutually skew as well as the lines  $B_j$; moreover each $A_i$ meets each $B_j$ except when $i=j$. 
    
    \noindent
     If the cubic surface $S$ is given as the image of the rational map
    $\xymatrix{
     \mu: \P^2 \ar@{-->}[r] & \P^3}$
     defined by the  linear system of plane cubic  curves 
   through six points $\{P_1, \dots, P_6\}$,  then a double-six  is implicitly defined by such a representation by letting $A_j$ be the transform of the point $P_j$ and $B_j$ the image of the plane conic $\theta_j$ containing $P_1, \dots, \widehat{P_j}, \dots, P_6$.  We have two morphisms:
   \[
   \xymatrix{
    \P^2& S \ar[l]_-{\pi_A} \ar[r]^-{\pi_B} & \P^2
   }
   \]
   $\pi_A$ is the contraction of  the lines $A_1, \dots, A_6$ and it is the inverse of $\mu$. Similarly $\pi_B$ contracts $B_1, \dots, B_6$ to points $R_1, \dots, R_6 \in \P^2$ and it is the inverse of the   rational map defined by the linear system of plane cubics  through $R_1, \dots, R_6$. 
    We have the following:
    
    \begin{theorem}\label{T:cre1}
    Each system of Cremona hexahedral equations of a nonsingular cubic surface $S$ defines a double-six  of lines on $S$. Conversely, the choice of a double-six of lines on $S$ defines a   system of Cremona hexahedral equations (\ref{E:cre1}) of $S$, which is uniquely determined up to replacing the coefficients $(\beta_0, \dots, \beta_5)$ by $(a+b\beta_0,\dots, a+b\beta_5)$ for some $a,b \in \k$, $b \ne 0$. 
    \end{theorem}

    We refer to \cite {D}, Theorem  9.4.6,  for the proof.   We need to point out   the following:
     
     \begin{corollary}\label{C:cre1}
     To a pair $(S,\Delta)$ consisting of a nonsingular cubic surface $S \subset \P^3$ and a  double-six of lines $\Delta$ on $S$, there is canonically associated a plane $\Xi \subset \P^3$   which is given by the equations 
          \begin{equation}\label{cre4}
     \cases{
    Z_0+Z_1+Z_2+Z_3+Z_4+Z_5 = 0  \cr
    \beta_0 Z_0+\beta_1 Z_1 + \beta_2 Z_2 + \beta_3 Z_3 + \beta_4 Z_4 + \beta_5 Z_5 = 0 \cr
    \beta_0^2Z_0+\beta_1^2Z_1 + \beta_2^2 Z_2 +\beta_3^2Z_3 +\beta_4^2Z_4+  \beta_5^2Z_5 = 0}
     \end{equation}
    \end{corollary}
    
    \proof By replacing in the third equation $\beta_i$ by $a+b \beta_i$.  with $b \ne 0$,  the plane $\Xi$ remains the same.  Therefore this plane depends only on the equations (\ref{E:cre1}) which in turn   depend only on $(S,\Delta)$.  \qed

    \begin{definition}\label{D:cocre}
    The plane $\Xi \subset \P^3$ will be called the \emph{Cremona plane} associated to the pair  $(S,\Delta)$.
    \end{definition}
    
   If a cubic surface $S \subset \P^3$ is given as the image of a linear system of plane cubic  curves 
   through six points $\{P_1, \dots, P_6\}$,  then a double-six    is implicitly selected by such a representation, and therefore $S$ can be given in $\P^5$ by   equations (\ref{E:cre1}).  In the first of the two papers \cite{Cob}  Coble found a parametrization of the  cubic surface  and of the constants $\beta_s$ depending on the points  $\{P_1, \dots, P_6\}$ so that   Cremona      hexahedral equations are satisfied. He defined six cubic polynomials $z_0,z_1, \dots, z_5 \in  H^0(\I_{\{P_1,\dots,P_6\}}(3))$ whose coefficients are multilinear in $P_1, \dots, P_6$, and such that:
   \begin{equation}\label{E:cob2}
 z_0+z_1+z_2+z_3+z_4+z_5 = 0 =   z_0^3+z_1^3+z_2^3+z_3^3+z_4^3+z_5^3
 \end{equation}
 identically.  
 In modern language the cubics $z_j$ span the 5-dimensional representation of $S_6$ which is called \emph{outer automorphism representation} (see  \cite{D}, \S 9.4.3).
 Moreover   Coble defined certain  multilinear polynomials $g_0, g_1, \dots, g_5$ in the $P_i$  satisfying the identity:
     $g_0+g_1+ \cdots + g_5 = 0$. 
    The representation of $S_6$ spanned by the $g_j$'s is the transpose of the outer automorphism representation, and it is obtained by tensoring with the sign representation.
    In fact Coble proves that the following identity holds:
    \begin{equation}\label{E:cob5}
    g_0z_0+g_1z_1+ \cdots + g_5z_5 = 0
    \end{equation}
     Putting together the identities (\ref{E:cob2}) and (\ref{E:cob5}) he then obtains the following:
    
    \begin{theorem}\label{C:cob1}
      The cubic surface $S\subset \P^3$ image of the rational map 
      $\xymatrix{\mu:\P^2 \ar@{-->}[r]&\P^3}$ defined by the linear system $|H^0(\I_{\{P_1,\dots,P_6\}}(3))|$   satisfies  the Cremona hexahedral equations (\ref{E:cre1}) with
      $\beta_s = g_s$.
      \end{theorem}

      Now consider: 
      \[
    C(P_1,\dots,P_6,X) = g_0^2z_0+g_1^2z_1+ \cdots + g_5^2z_5
    \]
    It is a  multihomogeneous polynomial  of degree 3 in $P_1,\dots,P_6,X$. The equation 
    $C(P_1,\dots,P_6,X) = 0$
    defines a cubic plane curve which is the pullback by $\mu$ of   the   Cremona plane $\Xi$   considered in Corollary \ref{C:cre1}.

     \begin{theorem}\label{T:morcob1}
     There is a constant $\lambda \ne 0$ such that the identity 
    \begin{equation}\label{E:morcob1}
     \Psi(P_1,\dots,P_6,P_7) = \lambda \ C(P_1,\dots,P_6,P_7)
     \end{equation}
    holds  for each 7-tuple of distinct points $P_1, \dots, P_6,P_7 \in \P^2$. 
      \end{theorem}
    
    \proof 
    Both $\Psi(P_1,\dots,P_6,P_7)$ and  $C(P_1,\dots,P_6,P_7)$ are cubic   SL$(3)$-invariants 
    of $P_1, \dots, P_6,P_7$ and skew-symmetric with respect to   the points.  
    Therefore it is enough to show that there is only one
skew-symmetric cubic  SL$(3)$-invariant of seven points, up to a constant factor.  This is proved in \cite{Cob}. We sketch a different approach to the proof.

\noindent
Let $S_n$ be the symmetric group of permutations on $n$ objects.
We denote by $\alpha=\alpha_1,\ldots,\alpha_k$
the Young diagram with $\alpha_i$ boxes in the i-th row and with $F^{\alpha}$
the corresponding representation of $S_n$.
Denote by $\Gamma^{\alpha}$ the Schur functor corresponding to the Young diagram $\alpha$.
The product group $SL(V)\times S_7$ acts in natural way on the vector space $S^3V\otimes\ldots\otimes S^3V$ (seven times).
We have the formula
$$S^3V\otimes\ldots\otimes S^3V=\oplus_{\alpha}\Gamma^{\alpha}(S^3V)\otimes F^{\alpha}$$
where we sum over all Young diagrams $\alpha$ with $7$ boxes, see \cite{Pro} ch. 9, Th. 3.1.4. 
The skew invariants are contained in the summand where
$\alpha=1^7$, indeed only in this case $F^{\alpha}$ is the sign representation of dimension one.
Correspondingly we can check, with a plethysm computation, that
$\Gamma^{1^7}(S^3V)=\wedge^7(S^3V)$ contains just one trivial summand of dimension one. This proves our result. \qed

Let us mention also that the same method gives another proof of the well known fact that every symmetric cubic invariant
of seven points is zero. 
      \begin{figure}
    \centerline{\includegraphics[width=70mm]{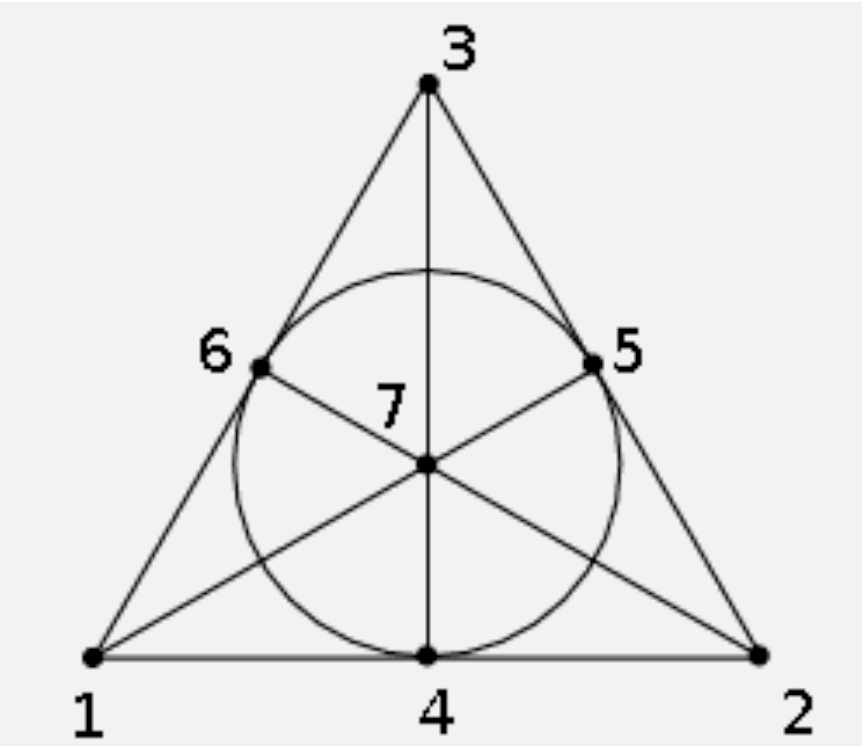}}
    \caption{The Fano plane}
    \end{figure}
      
    The symbolic expression of the Morley invariant is
$$
    |142||253||361||175||276||374||456|
$$    
Indeed skew-symmetrizing the previous expression over $S_7$  we get the Morley invariant $\Psi$. 
    The seven factors correspond to the seven lines of the Fano plane $\P^2_{\Z/2\Z}$, which is the smallest projective plane,  and  consists of seven points.
 Since the order of the automorphism group of the Fano plane is $168$, it is enough to consider just $7!/168=30$ summands, $15$ of them corresponding to even permutations, and the remaining $15$ corresponding to odd permutations. For another approach which uses Gopel functions see \cite{DO88}, ch. IX.

  \section{The   cubic of the seventh point}\label{S:cubic}

 The   analysis of \S \ref{S:morleinv} suggests the following. Given six distinct points $P_1, \dots, P_6 \in \P^2$, 
 no three of which are  on a line,   
 we can consider the  condition 
$\Psi(P_1,\dots,P_7) = 0$
 as defining a plane cubic curve $E_{P_1,\dots,P_6} \subset \P^2$ described by the seventh point $P_7$.  
    
    \begin{prop}\label{P:cubic1}
    In the above situation the cubic $E_{P_1,\dots,P_6}$ contains $P_1,\dots,P_6$.
    \end{prop}
    
    \proof  From the skew-symmetry of $\Psi(P_1, \dots, P_7)$  (Corollary \ref{C:pfaff1}) it follows that
    $\Psi(P_1, \dots, P_6,P_i)=0$ for all $i=1, \dots, 6$. \qed
    
    From this proposition it follows that  $E_{P_1,\dots,P_6}$   corresponds to a plane section 
     $\Xi \cap S$ of the cubic surface $S \subset \P^3$ determined by the linear system of cubics through $P_1,\dots,P_6$.  A description of the plane $\Xi$ will be given in \S \ref{S:plane}.  
     We will now look more closely at the curve $E_{P_1,\dots,P_6}$.

    \begin{prop}\label{P:cov1}
    Assume that $P_1, \dots , P_6$ are on a nonsingular conic $\theta$. Then $E_{P_1,\dots,P_6}$ is the closure of the locus of points $P_7\notin \theta$ such that the reducible conic $M^\seven_{P_7}$, where 
    $\seven = \{P_1, \dots , P_7\}$,    is conjugate to $\theta$. 
    \end{prop}

    \proof   By Remark \ref{R:mor1}  and by lower semicontinuity,   the condition    $P_7 \in E_{P_1,\dots,P_6}$   is equivalent to 
    $\dim(\langle M^\seven_{P_1}, \dots, M^\seven_{P_7} \rangle) \le 4$.
     By Proposition \ref{P:apolar0}  this is the condition  that the seven reducible conics $M^\seven_{P_i}$ are conjugate to $\theta$.  But when $i=1, \dots, 6$ this condition is automatically satisfied.   Therefore the only condition for $P_7 \in E_{a_1,\dots,a_6}$  is that 
 $M^\seven_{P_7}$  is conjugate to $\theta$, i.e. this condition defines  $E_{P_1,\dots,P_6}$.  \qed

 We have another description of  $E_{P_1,\dots,P_6}$, as follows.
 
 \begin{prop} \label{P:cov2} 
 Assume that $P_1, \dots , P_6$ are on a nonsingular conic $\theta$. Then $E_{P_1,\dots,P_6}$ is the cubic passing through  $P_1, \dots , P_6$ and apolar to $\theta$.
 \end{prop}

 \proof
 Note that the cubic $D$ passing through  $P_1, \dots , P_6$ and apolar to $\theta$ is unique by Proposition \ref{P:apol3}. Since both $D$ and $E_{P_1,\dots,P_6}$ are cubics it suffices to show that 
 $D \subset E_{P_1,\dots,P_6}$. By Prop. \ref{P:cov1} for this purpose it suffices to show that for each $P \in D$ and $P \ne P_1,\dots,P_6$, the cubic $G$ containing $P_1, \dots , P_6$ and singular at $P$ has the conic $M_P$ of principal tangents at $P$ conjugate to $\theta$.
 
 \noindent
 We may assume that $\theta = X_0^2 + 2X_1X_2$ and $P=(1,0,0)$.  Let 
 $L=aX_1+bX_2$ be any  line containing $P$.  From Proposition \ref{P:apol3} 
 it follows that   $\theta L$ is not apolar to $\theta$.  This means that  
 $D \notin \langle \theta X_1,\theta X_2 \rangle$ so that $\langle D, \theta X_1,\theta X_2 \rangle$ is the net of cubics through $P_1, \dots, P_6,P$. 
 
 \noindent
 If $D$ is singular at $P$ then $D=G$. In this case $M_P = \Delta_PD$  and this is conjugate to $\theta$ by Prop. \ref{P:apol2}.  Otherwise  
 \[
 D = \alpha_1(X_1,X_2)X_0^2 + \alpha_2(X_1,X_2)X_0 + \alpha_3(X_1,X_2) 
 \]
 with $\alpha_1 \ne 0$. Then
 \[
 G = D - \alpha_1\theta =  \alpha_2(X_1,X_2)X_0 + \alpha_3- 2\alpha_1X_1X_2
 \]
    so that $M_P = \alpha_2$.  On the other hand  
   $\Delta_PD = 2\alpha_1(X_1,X_2)X_0+ \alpha_2(X_1,X_2)$
    and the reducible conic joining $P$ to $ \Delta_PD \cap \Delta_P\theta$ is $\alpha_2$.  This is conjugate to $\theta$ by Prop. \ref{P:apol1}. \qed

    \begin{corollary}\label{C:cov1}
    The Morley invariant $\Psi(P_1,\dots,P_7)$ is irreducible and therefore the hypersurface ${\cal W} \subset (\P^2)^{(7)}$ is irreducible.
    \end{corollary}
    
    \proof 
    If $\Psi$ is reducible then the cubic $E_{P_1,\dots, P_6}$ is reducible for every    choice of 
    $P_1\dots P_6$.  But for a general choice of $P_1, \dots, P_6$ on a nonsingular conic $\theta$ the cubic $D$ passing through  $P_1, \dots , P_6$ and apolar to $\theta$ is irreducible by Proposition \ref{P:apol3} and 
    $D = E_{P_1,\dots, P_6}$ by the Proposition. \qed
    
     \begin{prop}\label{P:cov3}
      Assume   that   $P_1,\dots,P_6$ are not on a conic.  For each $i=1,\dots, 6$  
    let     $\theta_i$ be the conic containing all $P_j$'s except $P_i$ and let $D_i$ be the cubic containing   $P_1,\dots,P_6$ and apolar to $\theta_i$.  Denote by $Q_i \in \theta_i$ the sixth point of $D_i \cap \theta_i$.   Then      $E_{P_1,\dots,P_6}$ contains $Q_1, \dots, Q_6$. 
     \end{prop}

     \proof    Let $1 \le i \le 6$. Then  $D_i$ is apolar to $\theta_i$ and contains the six points 
     $P_1,\dots, \widehat P_i,\dots, P_6, Q_i$, which are on $\theta_i$. From Proposition \ref{P:cov2} it follows that  $D_i = E_{P_1,\dots, \widehat P_i,\dots, P_6, Q_i}$ and therefore
     \[
     \Psi(P_1,\dots,P_6,Q_i)=0,  \qquad i=1,\dots, 6
     \] 
     This means that  $E_{P_1,\dots,P_6}$ contains $Q_1, \dots, Q_6$.
      \qed

\section{A geometrical interpretation of the Cremona planes}\label{S:plane}

Consider a nonsingular cubic surface $S \subset \P^3$ and two skew lines $A,B \subset S$.   Denote by
   $ \xymatrix{f: A \ar[r] & B}$
   the double cover associating to $p \in A$ the point $f(p):= T_pS \cap B$ where $T_pS$ is the tangent plane to 
   $S$ at $p$.     
  Define
  $ \xymatrix{g: B \ar[r] & A}$
 similarly.  We   call $f$ and $g$ the \emph{involutory morphisms} relative to the pair of lines $A$ and $B$.
Let $p_1,p_2 \in A$  (resp. $q_1,q_2 \in B$) be the ramification points of $f$ (resp. $g$). 
Consider the pairs of branch points 
$f(p_1),f(p_2) \in B,  \  g(q_1), g(q_2) \in A$, 
and the new  morphisms
\[
\xymatrix{f': A \ar[r] & \P^1}, \qquad \xymatrix{g': B \ar[r] & \P^1}
\]
defined by the conditions that  $g(q_1), g(q_2)$ are ramification points of $f'$ and $f(p_1),f(p_2)$ are ramification points of $g'$.  
Let   $Q_1+Q_2$ (resp. $P_1+P_2$) be the   common divisor of  the two $g^1_2$'s on $A$ (resp. on $B$) defined by $f$ and $f'$ (resp. by  $g$ and $g'$). 
  The points 
 \[
 \bar P= g(P_1)=g(P_2) \in A, \qquad  \bar Q =f(Q_1)=f(Q_2)\in B
 \]
 are called  the \emph{involutory points} (relative to the pair of lines $A$ and $B$).

Consider six distinct points $P_1, \dots, P_6 \in \P^2$ not on a conic. Denote
by ${\cal A} = \{P_1, \dots, P_6\}$ and   by
$\xymatrix{\mu_{\cal A}: \P^2 \ar@{-->}[r]&   \P^3}$
 the   rational map defined by the linear system of cubics through ${\cal A}$.  On the   nonsingular cubic surface 
$S = \Im(\mu_{\cal A}) \subset  \P^3$
   consider the double-six of lines
   \[
   \Delta = (A_1, \dots, A_6; B_1, \dots, B_6)
   \]
   where  $A_1,\dots, A_6 \subset S$  are  the lines which are proper transforms under $\mu_{\cal A}$ of
 $P_1, \dots, P_6$ respectively, and $B_i   \subset S$ is the line which is the proper transform  of the
conic $\theta_i \subset \P^2$ containing $P_1, \dots, \widehat{P_i}, \dots, P_6$.  Consider the diagram
\[
   \xymatrix{
   \P^2& S \ar[l]_-{\pi_A} \ar[r]^-{\pi_B} & \P^2
    }
   \]
  where $\pi_A$ (resp. $\pi_B$) is the morphism which contracts the lines $A_1, \dots, A_6$  (resp. the lines $B_1, \dots, B_6$).  Let $R_i=\pi_B(B_i)\in \P^2$. Then $\pi_A$ is the inverse of $\mu_{\cal A}$  while $\pi_B$ is the inverse of
  $\xymatrix{\mu_{\cal B}: \P^2 \ar@{-->}[r]&   \P^3}$
where ${\cal B} = \{R_1, \dots, R_6\}$.
\noindent
Let $\bar P_i \in A_i, \  \bar Q_i \in B_i$ be the  involutory points  relative to the pair $A_i$ and $B_i$. We obtain twelve points
\[
\bar P_1, \dots, \bar P_6, \bar Q_1, \dots, \bar Q_6 \in S
\]
which are canonically associated to the double-six  $\Delta$.  

\begin{theorem}\label{T:plane1}
There is a unique plane  \  $\Xi \subset \P^3$ containing  the involutory points
\[
\bar P_1, \dots, \bar P_6, \bar Q_1, \dots, \bar Q_6 
\]
   Moreover $\Xi$ coincides with the Cremona plane (\ref{cre4}) associated to the pair $(S,\Delta)$ and 
   $\mu_{\cal A}^*(\Xi) = E_{P_1,\dots,P_6}$.
\end{theorem}

\proof
  We will keep the notations just introduced.   Fix $1 \le i \le 6$ and denote by 
  $\xymatrix{f_i: A_i \ar[r]&B_i}$ and by 
$\xymatrix{g_i: B_i \ar[r]&A_i}$ the  involutory morphisms. 
The  cubics $D$ through $P_1, \dots,  P_6$ and singular at $P_i$ form a pencil $\Lambda_i$ and are mapped by $\mu_{\cal A}$  on $S$ to
the  conics cut by the
planes containing $A_i$. 
Similarly,  the pencil $L_i$ of lines through $P_i$ is mapped by $\mu_{\cal A}$ on $S$ to the
  pencil  of  conics  cut by the planes
containing  $B_i$. 
It follows that $f_i$ can be interpreted in $\P^2$  as the map sending a line $\ell\in L_i$ to the sixth point of  $D \cap \theta_i$ where $D \in \Lambda_i$ is the cubic  having $\ell$ as a principal tangent. 
Therefore the ramification points of $f_i$ are the images  $p_1,p_2 \in A_i$ of the lines $\lambda_1,\lambda_2$   principal tangents of the two cuspidal cubics  $D_1,D_2 \in\Lambda_i$, and 
$f_i(p_1),    f_i(p_2)  \in \theta_i$
are the two sixth points of intersection of $D_1$, resp. $D_2$, with $\theta_i$.  

 \noindent
On the other hand  $g_i$  can be interpreted as associating to a point $q \in \theta_i$ the line $\langle P_i, q \rangle \in L_i$. 
The ramification points $q_1,q_2 \in \theta_i$  are  the  tangency points on the  two lines   
$\ell_1, \ell_2 \in L_i$  which are tangent   to $\theta_i$.

\noindent
It follows that  
$g_i(q_1)=\ell_1, \quad g_i(q_2) = \ell_2$.
Then clearly $\bar Q_i = \mu_{\cal A}(Q_i)\in B_i$,  where $Q_i \in \theta_i$ is the sixth point of $D \cap \theta_i$ and $D \in \Lambda_i$ is the cubic whose principal tangents are conjugate with respect to $\ell_1$ and $\ell_2$, or, equivalently, such that the reducible conic of its principal tangents is conjugate to $\theta_i$ (Prop. \ref{P:apol1}). 

\noindent
From Propositions \ref{P:cubic1} and  \ref{P:cov3} it follows that $E_{P_1,\dots,P_6}$ contains ${\cal A}$ and $Q_1, \dots, Q_6$, and clearly it is the only cubic curve with this property. Therefore there is a unique plane $\Xi \subset \P^3$ containing 
$\bar Q_1, \dots, \bar Q_6$ and $\mu_{\cal A}^*(\Xi)= E_{P_1,\dots,P_6}$.

\noindent
Reversing the role of $A_1, \dots, A_6$ and $B_1, \dots, B_6$ we can argue similarly using the rational map $\mu_{\cal B}$ instead of $\mu_{\cal A}$ to conclude that the points $\bar P_1, \dots, \bar P_6$ are contained in a unique plane $\Pi$ and that $\mu_{\cal B}^*(\Pi) = E_{R_1, \dots,R_6}$.  

\noindent
  From   Theorem \ref{T:morcob1} and the remarks  before it  we get that  both $\Xi$ and $\Pi$ coincide with the plane (\ref{cre4}) which is canonically associated to the double-six $\Delta$. In particular 
  $\Xi = \Pi$ and this concludes the proof.  \qed 
  
  \begin{remark}\label{R:plane}\rm
  The same proof as above shows that   the point $\bar P_i \in A_i$   corresponds to the line 
  $\tau_i := \langle P_i, z\rangle\in L_i$ joining $P_i$ with the pole $z$ with respect to $\theta_i$ of the line $\langle f_i(p_1), f_i(p_2) \rangle$.  
 Therefore   the theorem  implies  that the plane cubic curve $E_{P_1,\dots,P_6} \in H^0(\P^2,\I_\six(3))$   contains the points $Q_1, \dots, Q_6$ and its tangent   lines  at the points $P_1, \dots, P_6$ are $\tau_1, \dots, \tau_6$ respectively. 
\end{remark}

Since there are 36 double-six configurations of lines on a nonsingular cubic surface, we obtain 36 Cremona planes in $\P^3$ and, correspondingly, 36 cubic curves belonging to the linear system 
$|H^0(\P^2,\I_\six(3))|$. One of them is $E_{P_1,\dots,P_6}$.

        \begin{prop}\label{P:plane1}
        The 36   Cremona planes and, consequently,  the  corresponding 36 cubic curves in $H^0(\P^2,\I_\six(3))$, are pairwise distinct.
        \end{prop}
        
        \proof  We first remark that two double-sixes  always have a common line, we   call it $A_1$. The two corresponding skew lines $B_1$ and $B_1'$ (one for each double-six) are different. This follows from the explicit list of all the double-sixes, see  \cite{D}. It is well known that given a line $A_1$ on $S$, there are exactly $16$ lines on $S$ which are skew with $A_1$. Then it is enough to prove that the 16 involutory points on $A_1$ determined by the 16 lines skew with $A_1$ are all distinct. 
        Indeed, if two Cremona  planes corresponding to two different double-sixes coincide, then on their common line $A_1$ we get that  two points, among the 16 involutory points, should coincide.
        
We   computed explicitly these 16 points in a particular case, and we   got 16 distinct points, as we wanted. Let us sketch how this computation works.

\noindent
We start from a set ${\cal A}=\{P_1,\ldots, P_6\}$ of  six general points in the plane, and we consider the cubic surface $S =\Im(\mu_{\cal A})\subset \P^3$.
Denote by $A_i$ the exceptional divisor on   $S$ corresponding to  $P_i$.
The $16$  lines on $S$ which are skew with $A_1$ correspond to

a) the $5$ exceptional divisors $A_i$  for $i\ge 2$.

b) the $10$ lines joining two points among $P_2,\ldots, P_6$

c) the conic $\theta_1$ passing through $P_2,\ldots, P_6$.

\noindent
We have constructed the (five) involutory points on $A_1$ in case a), as explained in the proof of Theorem \ref{T:plane1}, starting from the points
$P_1=(1,0,0)$, $P_2=(0,1,0)$, $P_3=(0,0,1)$, $P_4=(1,1,1)$, $P_5=(2,3,5)$, $P_6=(11,13,29)$.
With the same points $P_i$, we have computed also the other ten points of case b) and the eleventh of case c). This last one corresponds to $\tau_1$ of Remark \ref{R:plane}.   
The resulting $16$ points on $A_1$ are all distinct.  \qed

  \section{The Geiser involution}\label{S:geiser}
    
    Our interest in configurations $\seven$ of seven distinct points in $\P^2$ comes from the classical Aronhold's construction of plane quartics  starting from such a $\seven$.   We refer to \cite{D} for details and proofs.  See also   
    \cite{EC29}, p.  319,  and  \cite{KW21}, p. 783, for  classical expositions, and \cite{jV08} for recent applications to vector bundles. 
    
    \noindent
    Let $\P$ be a projective plane and let $B \subset \P$ be a   nonsingular quartic.  Recall that an unordered 7-tuple $\T=\{t_1,\dots, t_7\}$ of bitangent lines of $B$ is called an \emph{Aronhold system}   if for all triples of distinct indices
    $1 \le i,j,k \le 7$ the six   points $(t_i \cup t_j \cup t_k) \cap B$ do not lie on a conic.
   Every nonsingular plane quartic has 288 distinct Aronhold systems of bitangents.

    \noindent
    Let  $\seven =\{P_1,\dots,P_7\} \subset \P(V)$ be seven distinct points such that no six of them are on a conic. The rational map 
    \[
    \xymatrix{
    \gamma_\seven: \P(V) \ar@{-->}[r] &|H^0(\I_\seven(3))|^\vee }
    \]
    defined by the net of cubics through $\seven$  is called the  \emph{Geiser involution} defined by $\seven$. It associates to a point $P \in \P(V)$ the pencil of cubics of the net containing $P$.  We have $\gamma_\seven(P) = \gamma_\seven(P')$ if and only if  $P$ and $P'$ are base points of the same pencil of cubics. Therefore $\gamma_\seven$ has degree 2 and is not defined precisely at the points of $\seven$.  
    We can identify the target of  $\gamma_\seven$  with $\P(V)^\vee$ by associating to $P$ the line  $\langle P,P'\rangle$ so that  $\gamma_\seven$ can be viewed as a rational map 
  $\xymatrix{
     \gamma_\seven: \P(V) \ar@{-->}[r] & \P(V)^\vee}$.
     Note that, conversely, from any general line $\ell \subset \P(V)$ we can recover $\{P,P'\}=\gamma_\seven^{-1}(\ell)$ as the  unique pair of point which are identified by the $g^2_3$ defined on $\ell$ by the net of cubics. 
     
     \noindent
          After choosing a basis of $V$   we obtain  a parametrization of the cubics of the net $|H^0(\I_\seven(3))|$ by associating to each $\xi \in \P(V)$ the cubic
     $S(\xi,X)$ (Lemma \ref{L:alg2}(ii)).  This defines an isomorphism:      
     \begin{equation}\label{E:isom1}
     \xymatrix{
     \P(V) \ar[r]^-\sim & |H^0(\I_\seven(3))|}
     \end{equation}
      The map   $\gamma_\seven$ can be described explicitly as follows.
     Let $P \in \P(V)$,  $P \notin \seven$.  Then
   $ \gamma_\seven(P)\in \P(V)^\vee$ is the line of $\P(V)$ given by the equation $S(\xi,P)=  0$ in coordinates $\xi$. This line parametrizes  the pencil of cubics of the net containing $P$ via the isomorphism (\ref{E:isom1}).

 \noindent
 Writing $S(\xi,X)= L_0(\xi)C_0(X)+L_1(\xi)C_1(X)+L_2(\xi)C_2(X)$, 
 the sextic
  \[
  \Sigma:  \left |{\partial C_j\over \partial X_h}\right | = 0
  \]
  is the  \emph{jacobian curve} of the net $|H^0(\I_\seven(3))|$, 
  i.e. the locus of double points of curves of the net.  
  
   All the properties of  $\gamma_\seven$ can be easily deduced by considering the Del Pezzo surface of degree 2 which is the blow-up of $\P(V)$ at $\seven$. The following proposition summarizes the main properties which we will need.

 \begin{prop}\label{P:geiser}
 \begin{itemize}
 \item[(i)]  The branch curve of the Geiser involution $\gamma_\seven$ is a   nonsingular quartic $B(\seven) \subset \P(V)^\vee$.

 \item[(ii)]
 The jacobian curve $\Sigma$ is the ramification curve of  $\gamma_\seven$.  
 
 \item[(iii)] The   points   $P_1,\dots, P_7$ are transformed into seven bitangent lines 
    $t_1, \dots, t_7$ of $B(\seven)$, which  form an Aronhold system.  The other 21 bitangent lines of $B(\seven)$ are the trasforms of the conics through five of the seven points of $\seven$. 
 \end{itemize}
 \end{prop}

     If we order the points $P_1, \dots,P_6,P_7$ we can consider the birational map
   $\xymatrix{
     \mu:    \P^2 \ar@{-->}[r] & S}$
      of  $\P^2$ onto a nonsingular cubic surface $S \subset \P^3$ defined by the linear system of cubics through $\{P_1, \dots, P_6\}$. The point $P_7$ is mapped to a point $O \in S$ and the ramification  sextic $\Sigma\subset \P^2$ is transformed into a sextic $\mu(\Sigma)\subset S$ of genus 3 with a double point at $O$. Then $\gamma_Z$ is the composition of $\mu$ with the projection of $S$ from $O$ onto the projective plane of lines through $O$. The quartic $B(Z)$ is the image of the sextic $\mu(\Sigma)$ under this projection.
    
   \separation \noindent
    Consider the following space:
    \[
    {\cal A}  := \left\{(B,T): \begin{array}{l}
    \hbox{$B\subset \P(V^\vee)$ is a n.s. quartic, and $T$ is an} \\
    \hbox{Aronhold system of bitangents of $B$}
    \end{array}\right\}
    \]
   It is   nonsingular and irreducible of dimension 14. We   have a commutative diagram of generically finite rational maps:
    \begin{equation}\label{E:aron1}
    \xymatrix{
     (\P^2)^{(7)}\ar@{-->}[r]^-{\cal T}\ar@{-->}[dr]^B & {\cal A}\ar[d]^\varphi \\
     &\P(S^4V)}
    \end{equation}
    where $B$ is the rational map associating to a   7-tuple $\seven$ of distinct points, no six of which are  on a conic, the quartic $B(\seven)$.  The map ${\cal T}$ associates to $\seven$ the pair
    ${\cal T}(\seven) = (B(\seven), \{t_1,\dots, t_7\})$
    where $t_i$ is the bitangent which is the image of $P_i$.  
     Since $B$ has finite fibres the image $B(\W) \subset \P(S^4V)$, where $\W\subset (\P^2)^{(7)}$ is defined in \ref{D:W},  is an open set of an irreducible hypersurface whose elements will be called \emph{Morley quartics.}

    \section{Morley quartics}\label{S:54}
    
    Consider  the irreducible hypersurface of Morley quartics
    \[
     \M:= \overline{B(\W)}\subset \P(S^4V)
    \]
    closure of the image of the hypersurface $\W \subset \P(V)^{(7)}$ under the rational map $B$.    Clearly  $\M$ is SL$(3)$-invariant.  In this section we will compute its degree.  
     
     \begin{theorem}\label{T:54}
     The hypersurface of Morley quartics $\M\subset \P(S^4V)$ has degree 54.
     \end{theorem}
      
\proof
Consider the projective bundle
$\xymatrix{
\pi: \P(\Q) \ar[r] &  \P^3}$
where $\Q= T_{\P^3}(-1)$ is the tautological quotient bundle which appears in the twisted Euler sequence: 
\[
\xymatrix{
0 \ar[r] &  {\cal O}(-1) \ar[r] &  {\cal O}^{\oplus 4} \ar[r] & \Q \ar[r]& 0 }
\]
For each $z \in \P^3$ the fibre $\pi^{-1}(z)$ is the projective plane of lines through $z$. Also consider the projective bundle
$
\xymatrix{
\beta: \P(S^4\Q^\vee) \ar[r] & \P^3}$.
  For each $z \in \P^3$ the fibre  $\beta^{-1}(z)$ is the linear system of quartics in  $\pi^{-1}(z)$.
The Picard group of $\P=\P(S^4\Q^\vee)$ is generated by $H= {\cal O}_{\P}(1)$ and by
  the pullback  $F$ of a  plane in $\P^3$.   
  Let 
$ \widetilde\M \subset \P(S^4\Q^\vee)$
be the $\beta$-relative hypersurface  of Morley quartics, and assume that it is given by a section of 
$aH+bF$.  Since it is invariant under the natural action of SL$(4)$ on $\P(S^4\Q^\vee)$, $\widetilde\M$ corresponds to a trivial summand of
  \[
 H^0(\P(S^4\Q^\vee),aH+bF)=H^0(\P^3,[S^a(S^4\Q)](b))
 \]
 which  exists if and only if $[S^a(S^4\Q)](b)$ contains $ {\cal O}$ as a summand. Since $\Q$ is homogeneous and indecomposable all the indecomposable summands of $[S^a(S^4\Q)](b)$ have the same slope. 
This is well known and can be easily deduced from the discussion at 5.2 of \cite{sR66}.
We get that
 $c_1\left([S^a(S^4\Q)](b)\right)=0$.
  Computing  the slope:
 \[
 \mu\left( [S^a(S^4\Q)](b) \right) = {4a \over 3}+b 
 \]
we deduce that 
$0 = 4a+3b$
and therefore $\widetilde\M= k(3H-4F)$ for some $k$. On the other hand  $\widetilde\M$ has relative degree $d$, where   $d$ is the degree of 
$\M\subset \P(S^4V)$.  Therefore $3k=d$ and it follows that $\widetilde\M$ has class $dH-{4d\over 3}F$.

\noindent
Let  $S \subset \P^3$ be a nonsingular cubic surface. To each $z\in S$   there is associated  the  quartic   branch curve of the rational projection 
$\xymatrix{S \ar@{-->}[r]& \pi^{-1}(z)}$ with center $z$.  This defines a section $s$ of $\beta$ over $S$:
\[
\xymatrix{\P(S^4\Q^\vee)_{|S} \ar[r]_-\beta & S \ar@/_1pc/[l]_s}
\] 
The pullback  $s^*\widetilde\M \subset S$ is the divisor of points $z\in S$ such that the branch curve of the projection of $S$ from $z$ is a Morley quartic. From    \ref{T:plane1} and \ref{P:plane1} it follows that  
$s^*\widetilde\M$ is a section of $ {\cal O}_S(\Xi_1 + \cdots + \Xi_{36})= {\cal O}_S(36)$, where $\Xi_1, \dots, \Xi_{36}$ are the  Cremona-Coble planes of $S$.   

Let's write an equation of $S$ as   $f(X,X,X)=0$, where  $f$ is a   symmetric trilinear form, and let $z\in S$.
The line through $z$ and $X$ is parametrized by $z+tX$ and meets $S$ where
$f(z+tX,z+tX,z+tX)$ vanishes.
We get:
\[
3tf(z,z,X)+3t^2f(z,X,X)+t^3f(X,X,X)=0
\]
which has a root  for $t=0$ and a residual double root when
\[
3f(z,X,X)^2-4f(z,z,X)f(X,X,X)=0
\]
which is a quartic cone with vertex in $z$.  From this expression we see that the section $s$ is quadratic in the coordinates of $z$.  
It follows that 
$s^*H =  {\cal O}_S(2)$. Therefore:
\[
 {\cal O}_S(36)=s^*\widetilde\M = s^*\left(dH-{4d\over 3}F\right) =  {\cal O}_S\left(2d-{4d\over 3}\right) =  {\cal O}_S(2d/3)
\]
 and we get $d=54$. \qed

 \begin{remark}\label{R:3/2}\rm
 The same   proof   shows that every invariant of a plane quartic
of degree $d$ gives a covariant of the cubic surface of degree $\frac{2d}{3}$. 
A classical reference for this statement is  \cite{Cob2}, p. 189. 
\end{remark}
 
 \section{Bateman configurations} \label{S:bateman}
    
    The configurations $\seven$ of seven points   in $\P(V)=\P^2$ no six of which are on a conic and for which the Morley invariant vanishes, i.e.  belonging to the irreducible hypersurface ${\cal W}$,   have a simple description which is due to   Bateman  \cite{hB14}.

    Consider  a nonsingular conic $\theta$, a cubic  $D$  and   the   matrix:
 \begin{equation}
\label{E:bat0}
A(\theta,D)=\pmatrix{\partial_0\theta&\partial_1\theta&\partial_2\theta \cr \partial_0D&\partial_1D&\partial_2D}
\end{equation}
 The   7-tuple of   points $\seven = \seven(\theta,D)$  defined by the maximal minors of this matrix   is called the \emph{Bateman configuration} defined by $\theta$ and $D$.    
   
   \noindent
   Note that  $\seven(\theta,D)$ consists of the points which have the same polar line with respect to $\theta$ and to $D$.

   \begin{lemma}\label{L:bateman1}
   Let $\seven = \seven(\theta,D)$ be the Bateman configuration defined by a nonsingular conic $\theta$ and by a cubic $D$.    If $D$ is general then   $\seven$ consists of seven distinct points  no six of which are on a conic.
 \end{lemma}

    \proof   We may assume that $\theta = X_0^2+X_1^2+X_2^2$. It suffices to prove the assertion in a special case.  If we take
    $D = X_0X_1X_2$ the maximal minors of $A(\theta,D)$ are: 
    \[
    C_0=X_0(X_1^2-X_2^2), \qquad C_1=X_1(X_2^2-X_0^2), \qquad C_2=X_2(X_0^2-X_1^2)
    \]
    and 
    \[
    \seven(\theta,D) = \{(1,0,0), (0,1,0), (0,0,1), (1,1,1), (1,-1,1), (1,1,-1), (-1,1,1)\}
    \]
    
    One easily checks directly that no six of the points of $\seven(\theta,D)$ are on a conic. This also follows from Lemma \ref{L:alg1}, because 
    $(\partial_0\theta,\partial_1\theta,\partial_2\theta) = (X_0,X_1,X_2)$ is the first row of $A(\theta,D)$ and has linearly independent entries. 
    \qed

       The Geiser involution associated to a general Bateman configuration $\seven=\seven(\theta,D)$ can be described as follows. Given  a point $P\in \P(V)$,  $P\notin \seven$,  then 
       $Q=(C_0(P),C_1(P),C_2(P))$ is the point of intersection of  the polar lines of $P$ with respect to $\theta$ and to $D$. The rational map 
       \[
       \xymatrix{\P(V) \ar@{-->}[r]& \P(V)}, \qquad \xymatrix{P\ar@{|->}[r]&Q}
       \]
      coincides with the Geiser involution $\gamma_\seven$ composed with  the identification  $\P(V^\vee)=\P(V)$ obtained thanks to the polarity with respect to $\theta$. In particular \emph{ the quartic $B(\seven(\theta,D))$ lies   in $\P(V)$.} 
        For each point $Q \in \P(V)$ we have
 $ \gamma_\seven^{-1}(Q) = \{P,P'\}$
 where   $P,P'$ are the points of intersection of the polar line of $Q$ with respect to $\theta$ with the polar conic of $Q$ with respect to $D$.

     We have:
   \begin{equation}\label{E:bat1}
    \begin{array}{cc}
   S(\xi,X) =  \left | \matrix{\partial_0\theta(\xi)&\partial_1\theta(\xi)&\partial_2\theta(\xi) \cr 
  \partial_0\theta(X)&\partial_1\theta(X)&\partial_2\theta(X) \cr \partial_0D(X)&\partial_1D(X)&\partial_2D(X)  }\right |   
   \end{array}
   \end{equation}
   and
   \[
   \begin{array}{cc}
   M(\xi,X) =  \left | \matrix{\partial_0\theta(\xi)&\partial_1\theta(\xi)&\partial_2\theta(\xi) \cr 
  \partial_0\theta(X)&\partial_1\theta(X)&\partial_2\theta(X) \cr 
  \Delta_\xi\partial_0D(X)&\Delta_\xi\partial_1D(X)&\Delta_\xi\partial_2D(X)}\right |   
   \end{array}
     \]
   
   \begin{theorem}\label{T:bateman}
   Let $\seven(\theta,D)$ be a Bateman configuration. Then, for every $\xi \in \P^2$, the conic $M(\xi,X)$ is conjugate to   $\theta$.  
   \end{theorem}
    
    \proof 
    We can change coordinates and assume that $\theta=X_0^2 + 2X_1X_2$,
  so that its dual is 
  $ \theta^* =  \partial_0^2+2\partial_1\partial_2$,
   and we must show that 
  \begin{equation}\label{E:diffid2}
   P_{\theta^*}(M(\xi,X))=0
   \end{equation}
   identically, where 
   \[
   \begin{array}{cc}
   M(\xi,X) =  \left | \matrix{ \xi_0 &  \xi_2& \xi_1 \cr 
  X_0  & X_2   &  X_1  \cr 
  \Delta_\xi\partial_0D(X)&\Delta_\xi\partial_1D(X)&\Delta_\xi\partial_2D(X)  }\right |   = \\ \\
   \xi_0(X_2\Delta_\xi \partial_2D - X_1\Delta_\xi \partial_1D) - \xi_2(X_0\Delta_\xi \partial_2D - X_1\Delta_\xi \partial_0D) +
  \xi_1(X_0\Delta_\xi \partial_1 D-X_2\Delta_\xi\partial_0 D)
  \end{array}
  \]
  Writing the cubic polynomial defining $D$   as:
   \[
   \sum_{0 \le i\le j\le k\le 2} \beta_{ijk}X_iX_jX_k
   \]
   an easy computation shows that
  \[
  \begin{array}{ll}
  M(\xi,X) = &(X_1X_2-X_0^2)(\beta_{002}\xi_0\xi_2 - \beta_{001}\xi_0\xi_1+\beta_{022}\xi_2^2 - \beta_{011}\xi_1^2) \\ \\
  & + \quad (\hbox{terms not involving $X_0^2$ and $X_1X_2$})
  \end{array}
  \]
  and (\ref{E:diffid2}) follows immediately.  \qed 
  
  \begin{corollary}\label{C:bateman1}
  If $\seven(\theta,D) = \{P_1,\dots,P_7\}$  is a Bateman configuration  of distinct points no six of which are on a conic,    then $\Psi(P_1,\dots,P_7)=0$.
   In other words, the image of the rational map
    \[
    \xymatrix{
    \seven: \P(S^2V^\vee) \times \P(S^3V^\vee) \ar@{-->}[r] & \P(V)^{(7)}=(\P^2)^{(7)}}\]
    which associates to a general pair $(\theta,D)$ the Bateman configuration $\seven(\theta,D)$, is contained in ${\cal W} \subset  (\P^2)^{(7)}$ (see Definition \ref{D:W}).
  \end{corollary}
  
  \proof From the theorem it follows that all the conics $M(\xi,X)$, $\xi\in \P^2$,  are contained in the hyperplane of conics conjugate to the   conic $\theta$. This implies that the skew-symmetric form
   \[
   M: S^2V^\vee \times S^2V^\vee \to \k
   \]
   is degenerate, hence its pfaffian vanishes. But since no six of the points of $\seven(\theta,D)$ are on a conic, we have
   ${\cal Q}(P_1,\dots,\widehat{P}_i,\dots, P_7) \ne 0$ for all $i=1,\dots, 7$. Then the conclusion  follows from the factorization 
   (\ref{E:pfaff1}). \qed
  
  \begin{definition}\label{D:diffid}
  Identity (\ref{E:diffid2})  is called \emph{Morley's differential identity.}
  \end{definition}

  Corollary \ref{C:bateman1} shows in particular that Bateman configurations are not the most general 7-tuples of points because they are in  ${\cal W}$. The corollary does not exclude that Im$(\seven)$, i.e. the locus of Bateman configurations, is contained in a proper closed subset of ${\cal W}$.  We will show in \S \ref{S:luroth} that  the  Bateman configurations actually fill a dense open subset of   ${\cal W}$, i.e. they depend on 13 parameters and not less.

  \section{L\"uroth quartics} \label{S:luroth}
  
  A configuration consisting of five lines in $\P(V)$, three by three linearly independent,  
    together with the ten double  points of their  union will be called a \emph{complete pentalateral}.  The ten nodes of their union are called \emph{vertices} of the complete pentalateral.

       \begin{definition}\label{D:lur}
  A \emph{L\"uroth quartic} is a nonsingular   quartic  $B \subset \P(V)$ which is circumscribed to a complete pentalateral, i.e. which contains its  ten vertices.
   \end{definition}

   Consider the incidence relation $\widetilde{\L} \subset  \P(S^4V^\vee)\times \P(V^\vee)^{(5)}$
   described as follows:
   \[
   \widetilde{\L} := \left\{(B,\{\ell_0,\dots,\ell_4\}):
   \begin{array}{l}
   \hbox{$\{\ell_0,\dots,\ell_4\}$ is a complete pentalateral and} \\ 
   \hbox{$B$ is a n.s. quartic circumscribed to it}\end{array}\right\}   
   \]
   and the projections:
   \[
   \xymatrix{
    \P(S^4V^\vee)&  \widetilde{\L} \ar[l]_-{q_1}\ar[r]^-{q_2}&\P(V^\vee)^{(5)} 
  }
   \]
   Clearly   $q_1(\widetilde{\L}) \subset \P(S^4V^\vee)$ is the locus of L\"uroth quartics. The following facts are well known (see \cite{gO07}):
   \begin{itemize}
   \item[(i)] $q_2$ is dominant with general fibre of dimension 4 and $ \widetilde{\L}$ is   irreducible of dimension 14.
   
   \item[(ii)] The general fibre of $q_1$ has dimension 1. This means that every L\"uroth quartic has infinitely many inscribed pentalaterals. Moreover 
   \[
   \L := \overline{q_1(\widetilde{\L})}  \subset \P(S^4V^\vee)
   \]
   is an   $SL(3)$-invariant  irreducible hypersurface which is called     \emph{the L\"uroth hypersurface}.  
   \end{itemize}
      Given  a complete pentalateral $\{\ell_0,\dots,\ell_4\} \in \P(V^\vee)^{(5)}$  , the   fibre $q_2^{-1}(\{\ell_0,\dots,\ell_4\})$ is the linear system of quartics circumscribed to it.  It consists of the quartics of the form:
      \begin{equation}\label{E:lur1}
     \sum_{k=0}^4 \lambda_k\ell_0\cdots \hat\ell_k \cdots \ell_4 =0
     \end{equation}
     as  $(\lambda_0,\dots,\lambda_4) \in \P^4$. Another way  of describing  a general element of $q_2^{-1}(\{\ell_0,\dots,\ell_4\})$ is under the form:
      \begin{equation}\label{E:lur2}
   \sum_{k=0}^4 {1 \over l_k \ell_k}  =0
  \end{equation}
  as $(l_0,\dots,l_4) \in \P^4$.  The two descriptions are of course related by a Cremona transformation of $\P^4$. 
  We have moreover the following elementary property:
  \begin{itemize}
   \item[(iii)] Any given $(B,\{\ell_0,\dots,\ell_4\}) \in \widetilde{\L}$ is uniquely determined by any of the five pairs
   \[
   (B,\{\ell_0,\dots,\hat\ell_k, \dots,\ell_4\}) \in\P(S^4V^\vee)\times \P(V^\vee)^{(4)}
   \]
    consisting of the quartic $B$ and of four of the five lines of the pentalateral.
    \end{itemize}

  We will   need the following result,   due to R. A. Roberts.
 
\begin{theorem}[\cite{rR89}]\label{T:roberts}
  Let $\theta$ be a nonsingular conic  and $D$ a general cubic. Then there are   lines $\ell_1,\ell_2,\ell_3,\ell_4$,  three by three linearly independent,  and constants 
  $a_1,a_2,a_3,a_4, b_1,b_2, b_3, b_4$ such that
  \begin{equation}\label{E:rob1}
  \begin{array}{ll}
   \theta = a_1\ell_1^2+a_2\ell_2^2+a_3\ell_3^2+a_4\ell_4^2 \\ 
  D = b_1 \ell_1^3 + b_2\ell_2^3 + b_3\ell_3^3 +b_4\ell_4^3  
 \end{array}
  \end{equation}
  The four lines are   uniquely determined and each of the two 4-tuples of constants is   uniquely determined up to a constant factor.
  \end{theorem}
  
  \proof
  The line-conics which are simultaneously apolar to $\theta$ and to $D$ form at least a pencil, because 
    being apolar to $\theta$, resp, to $D$, is one condition, resp. three conditions, for a line-conic. Moreover, for a general choice of $(\theta,D)$, these conditions are independent. In fact, taking $(\theta,D)$ as in (\ref{E:rob1}), and letting $\Sigma$ be a line-conic belonging to the pencil with base the lines $\ell_1,\ell_2,\ell_3,\ell_4$, we have:
    \[
   P_\Sigma(\theta) = 2a_1 \Sigma(\ell_1) + 2a_2 \Sigma(\ell_2) + 2a_3 \Sigma(\ell_3) + 2a_4 \Sigma(\ell_4) = 0
   \]
    Similarly  $P_\Sigma(D)=0$. 
   In other words $\Sigma$ is apolar to both $\theta$ and $D$. On the other hand, it is clear that there are no other  line-conics apolar to $\theta$ and $D$. Now the theorem follows   from Proposition 4.3 of 
  \cite{DK93}. \qed

     Theorem \ref{T:roberts} can be conveniently rephrased as follows:
     
     \begin{corollary}\label{C:roberts}
     By associating to a pair $(\theta,D)\in \P(S^2V^\vee) \times \P(S^3V^\vee)$ consisting of a nonsingular conic and a general cubic the data
     \[
     ((\ell_1,\ell_2,\ell_3,\ell_4),(a_1,a_2,a_3,a_4),( b_1,b_2, b_3, b_4))
     \]
    given by Theorem \ref{T:roberts} one obtains a birational map:
    \[
    \xymatrix{
    R: \P(S^2V^\vee) \times \P(S^3V^\vee) \ar@{-->}[r] & (\P^{2\vee})^4 \times \P^3 \times \P^3 }
    \]
    \end{corollary}

  We have the following remarkable   result.

 \begin{theorem}[Bateman \cite{hB14}]\label{T:bateman2} 
 Consider a nonsingular conic  $\theta$ and a general  cubic $D$,  and represent them as
     \[
     \begin{array}{ll}
   \theta =  a_1\ell_1^2+a_2\ell_2^2+a_3\ell_3^2+a_4\ell_4^2, & 
  D =  b_1\ell_1^3+b_2\ell_2^3+b_3\ell_3^3  + b_4 \ell_4^3
 \end{array}
 \]
 according to Theorem \ref{T:roberts}. Then the plane quartic $B=B(\seven(\theta,D))$ associated to the Bateman configuration $\seven(\theta,D)$ is a L\"uroth quartic and $\ell_1,\ell_2,\ell_3,\ell_4$ are four lines of a complete pentalateral inscribed in $B$.
 \end{theorem}

 \proof
  Let $X_1,\dots, X_4$  be homogeneous coordinates   in $\P^3$. We may identify  $\P(V)$ with the plane 
 $H\subset \P^3$ of equation    $\sum_i X_i=0$. After a change of coordinates in $H$ we may further assume that  $\ell_1,\ell_2,\ell_3,\ell_4$ are respectively the lines  $X_i=0$,  $i=1,\dots, 4$.
  With this convention we have:
 \[
 \theta= \sum_{i=1}^4 a_iX_i^2, \qquad  D = \sum_{i=1}^4b_iX_i^3 
 \]
 and we may assume that the constants  $a_i,b_i $ are   all non-zero.  Let $Q=(y_1,\dots,y_4) \in H$ be a (variable) point. The polar line of $Q$ w.r. to $\theta$ is:
 \begin{equation}\label{E:3}
  \sum_k \ a_ky_k X_k
 \end{equation}
 Similarly the polar conic of $Q$ w.r. to $D$ is:
  \begin{equation}\label{E:4}
  \sum_k \ b_ky_k X_k^2 
 \end{equation}
 Assume that the line (\ref{E:3}) is tangent to the conic (\ref{E:4}) at the point $P=(z_1,\dots,z_4)$. Then its equation must be equivalent to the equation
 $\sum_k  b_ky_kz_kX_k =0$. 
 This means that there are constants $(\lambda,\mu) \ne (0,0)$ such that
 \[
   \sum_k  b_ky_kz_kX_k  = \lambda\left[\sum_k \ a_ky_k X_k\right]+ \mu\left[\sum_k X_k\right]
 \]
 or, equivalently:
 \[
 b_ky_kz_k = \lambda a_ky_k  + \mu,  \qquad k=1,2,3,4.
 \]
 Since $P\in H$ we find:
 \[
 0 = \sum_kz_k = \lambda\left[\sum_k{a_k\over b_k}\right] + \mu\left[\sum_k{1\over b_ky_k}\right]
 \]
 Using the fact that $P$ belongs to the polar line (\ref{E:3}) we also deduce that:
 \[
 0 = \sum_k a_ky_k z_k = \lambda\left[\sum_k{a_k^2y_k\over b_k}\right] + \mu\left[\sum_k{a_k\over b_k}\right]
 \]
 These two identities imply that:
 \[
 \vrule\begin{array}{cc}
 \sum{a_k\over b_k} & \sum{1\over b_ky_k} \\ \\
 \sum{a_k^2y_k\over b_k} & \sum{a_k\over b_k}
 \end{array}\vrule \ =0
 \]
 or, equivalently:
 \begin{equation}\label{E:5}
\left( \sum_{k=1}^4{a_k\over b_k}\right)^2 = 
\left(\sum_{k=1}^4{a_k^2y_k\over b_k}\right)\left( \sum_{k=1}^4{1\over b_ky_k}\right) 
 \end{equation}
 Now let's define 
 \begin{equation}\label{E:lur3}
 L =  -\left( \sum_{k=1}^4{a_k\over b_k}\right)^{-2}\left(\sum_{k=1}^4{a_k^2y_k\over b_k}\right)
 \end{equation}
 Then $L$ is a linear form in the coordinates $y_1,\dots, y_4$ of $Q$, and the identity (\ref{E:5}) is equivalent to:
 \begin{equation}\label{E:6}
 \sum_{k=1}^4 {1\over b_ky_k}+ {1\over L} =0
 \end{equation}
 This is the equation of a L\"uroth quartic in the coordinates of $Q$. 
 \qed

 \begin{corollary}\label{T:bateman3}
  There is a dominant, generically finite, rational map
  \[
 \xymatrix{
  \widetilde{\seven}:\P(S^2V^\vee) \times \P(S^3V^\vee) \ar@{-->}[r] &\widetilde{\L} }
  \]
   such that $q_1(\widetilde{\seven}(\theta,D)) = B(\seven(\theta,D))$.
     In particular:
  \begin{itemize}
  
  \item[(i)]  The general L\"uroth quartic is of the form $B(\seven(\theta,D))$ for some 
  $(\theta,D)$.
  
  \item[(ii)]  The rational map 
  $\xymatrix{\seven: \P(S^2V^\vee) \times \P(S^3V^\vee)\ar@{-->} [r]& {\cal W}}$
   of Corollary \ref{C:bateman1} is dominant.
  
  \end{itemize}
   \end{corollary}

     \proof 
      Consider  a pair $(\theta,D)$ consisting of a nonsingular conic and a general cubic. By Theorem \ref{T:bateman2}  the quartic $B(\seven(\theta,D))$ is L\"uroth. Moreover, again by Theorem \ref{T:bateman2},  the lines $\ell_1, \dots,\ell_4$ associated to $(\theta,D)$ by Theorem \ref{T:roberts} are components of a complete pentalateral inscribed in $B(\seven(\theta,D))$. Then the map $\widetilde{\seven}$ is defined by associating to $(\theta,D)$ the pair
  $(B(\seven(\theta,D)), \{\ell_0,\ell_1, \dots,\ell_4\}) \in \widetilde{\L}$,
  where $\ell_0$ is the fifth line of the complete pentalateral inscribed in $B(\seven(\theta,D))$ having     $\ell_1, \dots,\ell_4$ as components.  
  
  \noindent
   Consider a general  $(\theta,D) \in \P(S^2V^\vee) \times \P(S^3V^\vee)$ and let
  \[
  (B, \{\ell_0,  \dots, \ell_4\}) = \widetilde{\seven}(\theta,D)
  \]
  By the definition of $\widetilde{\seven}$ it follows that for some $0 \le k \le 4$ the lines
  $\ell_0, \dots, \hat\ell_k,\dots,\ell_4$ are simultaneously apolar to $\theta$ and to $D$.  We may assume that $k=0$ and choose coordinates so that  $\ell_1+\ell_2+\ell_3+\ell_4 = 0$.  Then from the proof of Theorem   \ref{T:bateman2} it follows that
 $B$ has equation of the form (\ref{E:6}), which in our notation takes the form:
 \[
 {1\over b_1\ell_1}+{1\over b_2\ell_2}+{1\over b_3\ell_3} + {1 \over b_4\ell_4} + {1\over L} = 0
 \]
 where $b_1, \dots, b_4$ are the uniquely defined non-zero coefficients such that 
 $D = b_1\ell_1^3+ \cdots + b_4\ell_4^3$
and $L$ is a linear combination of $\ell_1, \dots, \ell_4$  given by   (\ref{E:lur3}), which now takes the form:
  \begin{equation}\label{E:lur5}
  L = -\left( \sum_{k=1}^4{a_k\over b_k}\right)^{-2}\left(\sum_{k=1}^4{a_k^2\ell_k\over b_k}\right)
\end{equation}
where $a_1, \dots, a_4$ are the uniquely determined non-zero coefficients such that
$\theta = a_1\ell_1^2 + \cdots + a_4\ell_4^2$.   
 From these expressions it follows that $D$ is uniquely determined by  $(B, \{\ell_0,  \dots, \ell_4\})$.
 The coefficients of the linear combination (\ref{E:lur5}) are rational functions of $a_1,\dots,a_4$, homogeneous of degree zero which can be interpreted as follows. Let
 $\xymatrix{
 \P^3 \ar@{-->}[r] & \P^4}$
 be defined by sending:
 \[ 
 \xymatrix{
 (a_1,\dots,a_4) \ar@{|->}[r] &
 \left(\left(\sum_{k=1}^4{a_k\over b_k}\right)^2,{a_1^2\over b_1}, \dots,{a_4^2\over b_4} \right)}
 \]
  This being the composition of a Veronese map with a projection, is finite on its set of definition. From this remark it follows that, given $b_1, \dots, b_4$, there are finitely many $a_1,\dots, a_4$ defining $L$. This shows that $(\theta,D)$ is isolated in 
  $ \widetilde{\seven}^{-1}(B, \{\ell_0,  \dots, \ell_4\})$. 
  Since  its domain and range    are irreducible of dimension  14, this proves that    $\widetilde{\seven}$ is dominant and  generically finite.  
  The assertions (i), (ii) are now obvious, recalling that ${\cal W}$ is irreducible (see Corollary \ref{C:cov1}).
  \qed

Our final and main result  is:

\begin{theorem}\label{T:main}
The hypersurface $\L \subset \P(S^4V^\vee)$ has degree 54.  
\end{theorem}

\proof
From Corollary \ref{T:bateman3} we deduce that the hypersurface of L\"uroth quartics can be identified  with the hypersurface 
$\M$ of Morley quartics. In particular
$\deg(\L) = \deg(\M)=54$.
\qed

\begin{remark}\rm\label{R:pentath}
We recall the following construction from \cite{D}. Given an Aronhold system $\{t_1,\ldots ,t_7\}$ of bitangents for a nonsingular plane quartic $B$, consider them as odd theta characteristic, and let $H$ be the divisor on $C$ cut by a line. The $35$ divisors $t_i+t_k+t_k-H$, for $1\le i<j<k\le 7$, define $35$ distinct  even theta characteristic. Since there are $36$ even theta characteristic, we may denote the remaining one by $t(t_1,\ldots ,t_7)$.

\noindent
We come back to the Cremona  planes of   \S   \ref{S:plane}. Given a nonsingular cubic surface $S \subset \P^3$ and a double-six   $\Delta =(A_1,\ldots, A_6; B_1,\ldots, B_6)$,   a point $P\in S$  belongs to the corresponding Cremona plane $\Xi$ if and only if the quartic branch curve $B$ of the rational projection with center $P$ is L\"uroth and the vertices $v_1, \dots, v_{10} \in B$ of an inscribed pentalateral  satisfy
$v_1+ \cdots+v_{10} \in |2H+t(t_1,\ldots ,t_7)|$,
   where 
$\{t_1,\dots,t_7\}$ is the Aronhold system consisting of the bitangents which are projections of $A_1, \dots, A_6,P$. 
It is natural to call  $t(t_1,\ldots ,t_7)$   \emph{the pentalateral even  theta-characteristic} on $B$.

\end{remark}

 \noindent
  \textsc{g. ottaviani} -
  Dipartimento di Matematica ``U. Dini'', Universit\`a di Firenze, viale Morgagni 67/A, 50134 Firenze (Italy). e-mail: \texttt{ottavian@math.unifi.it}

  \separation
  \noindent
\textsc{e. sernesi} - 
 Dipartimento di Matematica,
 Universit\`a Roma Tre, 
 Largo S.L. Murialdo 1-
 00146 Roma (Italy).  e-mail: \texttt{sernesi@mat.uniroma3.it}


\begin{thebibliography}{999}
  
 \bibitem{Bar} W. Barth:  Moduli of vector bundles on the projective plane , \emph{Invent. Math.} {\bf 42} (1977), 63--92.

 \bibitem{hB14}
  H. Bateman: The Quartic Curve and its Inscribed Configurations,  \emph{American J. of Math.} 36 (1914), 357-386. 
  
  \bibitem{Br}
  A. Brouwer: The invariant theory of plane quartics, avail. at 
  \begin{verbatim}http://www.win.tue.nl/~aeb/math/ternary_quartic.html \end{verbatim}
 

  \bibitem{CGO}
  C. Ciliberto, A.V. Geramita, F. Orecchia: Remarks on a theorem of Hilbert-Burch, in 
  \emph{The Curves Seminar  at Queen's, IV}, Queen's Papers in Pure and Applied Math.  76 (1986). 
  
  \bibitem{Cob} A. B. Coble:  Point sets and allied Cremona groups, I, II \emph{Trans. AMS,}  16 (1915),   17 (1916).
  
  \bibitem{Cob2} A. B. Coble: \emph{Algebraic Geometry and Theta Functions} (reprint of the 1929 edition), AMS Colloquium Publications v. 10 (1982).
  
  
   \bibitem{Cre} L. Cremona: Ueber die Polar-Hexaeder bei den Fl\"achen dritter Ordnung, \emph{Math. Annalen} 13 (1877)
 
\bibitem{D} I. Dolgachev: \emph{Topics in Classical Algebraic Geometry,}   
avail. at 
 \begin{verbatim}http://www.math.lsa.umich.edu/~idolga/lecturenotes.html. \end{verbatim}

\bibitem{Do2} 
I. Dolgachev: \emph{Lectures on invariant theory.} London Mathematical Society Lecture Note Series, 296. CUP, Cambridge, 2003.

\bibitem{DK93}
I. Dolgachev, V. Kanev: Polar covariants of plane cubics and quartics, \emph{Advances in Math.} 98 (1993), 216-301.
  
   \bibitem{DO88} I. Dolgachev, D. Ortland: \emph{Point Sets in Projective Spaces and Theta Functions},  Asterisque 165 (1988).
  
  \bibitem{dE05}D. Eisenbud: \emph{The Geometry of Syzygies}, Springer GTM v. 229 (2005).

\bibitem{EPS} G. Ellingsrud, J. Le Potier,   S.A. Stromme: Some Donaldson invariants of ${\bf C}\P^2$,  in \emph{Moduli of vector bundles (Sanda, 1994; Kyoto, 1994)},  33--38, Lecture Notes in Pure and Appl. Math., 179, Dekker, New York, 1996.

\bibitem{EC29} F. Enriques, O. Chisini: \emph{Teoria geometrica delle equazioni e delle funzioni algebriche},  vol. I, Zanichelli, Bologna 1929.
 
\bibitem{KW21} A. Krazer, W. Wirtinger: Abelsche Funktionen und allgemeine Thetafunktionen. In \emph{Encyklopadie der mathematischen Wissenschaften} II (2), Heft 7, 604-873,  Leipzig 1921. 

\bibitem{PT} J. Le Potier, A. Tikhomirov:  Sur le morphisme de Barth,
\emph{Ann. Sci. \'Ecole Norm. Sup.} (4) {\bf 34} (2001), no. 4, 573--629.

 

\bibitem{LQ} W.P. Li, Z. Qin:  Lower-degree Donaldson polynomial invariants of rational surfaces.  \emph{J. Algebraic Geom.}  2  (1993),  no. 3, 413--442.

\bibitem{jL68}
J. L\"uroth: Einige Eigenshaften einer gewissen Gathung von Curves vierten Ordnung, 
\emph{Math. Annalen} 1 (1868), 38-53.
  
  \bibitem{fM14}
  F. Morley: On the L\"uroth Quartic Curve,   \emph{American J. of Math.} 41 (1919), 279-282. 
  
  
  \bibitem{gO07}
 G. Ottaviani: Symplectic bundles on the plane, secant varieties and L\"uroth quartics revisited, Math.AG/0702151,
  in Quaderni di Matematica, vol. 21 (eds. G. Casnati, F. Catanese, R. Notari), \emph{Vector bundles and Low Codimensional Subvarieties:
State of the Art and Recent Developments,} Aracne, 2008.
 
\bibitem{Pro} 
 C. Procesi:  \emph{Lie Groups,} Springer New York, 2007. 
  
\bibitem{sR66}
S. Ramanan: Holomorphic vector bundles on homogeneous spaces, \emph{Topology} 5 (1966), 159-177.
  
  
  \bibitem{rR89}
  R. A. Roberts: Note on the Plane Cubic and a Conic, \emph{Proc. London Math. Soc.} (1) 21 (1889),62-69. 

  
  \bibitem{jV08}
 J. Vall\`es: Fibr\'es vectoriels de rang deux sur $\P^2$ provenant d'un revetement double, arXiv:0806.3355, to appear in Ann.  Inst. Fourier.

 
  \end{thebibliography}
\end{document}